\documentclass[preprints,article,accept,moreauthors]{my-mdpi} 

\firstpage{1} 
\pubvolume{11}
\issuenum{6}
\articlenumber{1511}
\pubyear{2023}
\copyrightyear{2023}
\externaleditor{Academic Editor: Ellina Grigorieva}
\datereceived{31 January 2023} 
\daterevised{13 March 2023} 
\dateaccepted{16 March 2023} 
\datepublished{20 March 2023} 
\hreflink{https://doi.org/\linebreak 10.3390/math11061511} 
\doinum{10.3390/math11061511}

\usepackage{amsfonts,caption}
\usepackage{bbm}
\usepackage{booktabs,multirow}
\usepackage{appendix}
\usepackage[labelformat=simple]{subcaption}

\DeclareCaptionLabelFormat{subcaptionlabel}{\normalfont(\textbf{#2}\normalfont)}
\usepackage{alltt,xcolor}

\newcommand{\Dl}{{}^{\scriptscriptstyle C}_{\scriptscriptstyle 0}\!D^{\alpha}_{\scriptscriptstyle t}}
\newcommand{\Dll}{{}^{\scriptscriptstyle C}_{\scriptscriptstyle 0}\!D^{\alpha}_{\scriptscriptstyle t'}}
\newcommand{\DlRL}{{}_{\scriptscriptstyle 0}D^{\alpha}_{\scriptscriptstyle t'}}
\newcommand{\DrRL}{{}_{\scriptscriptstyle t}D^{\alpha}_{\scriptscriptstyle t_f}}
\newcommand{\Dlmm}{{}^{\scriptscriptstyle C}_{\scriptscriptstyle 0}\!D^{-\alpha}_{\scriptscriptstyle t}}

\Title{Numerical Fractional Optimal Control of Respiratory Syncytial Virus Infection in Octave/MATLAB}

\TitleCitation{Numerical Fractional Optimal Control of Respiratory Syncytial Virus Infection in Octave/MATLAB}

\Author{Silv\'erio Rosa $^{1,\dagger}$\orcidA{}
and Delfim F. M. Torres $^{2,}$*${^{,\dagger}}$\orcidB{}}

\AuthorCitation{Rosa, S.; Torres, D.F.M.} 

\address{$^{1}$ \quad Instituto de Telecomunica\c{c}\~{o}es and Department of Mathematics, 
University of Beira Interior, \linebreak 
6201-001 Covilh\~a, Portugal; rosa@ubi.pt \\
$^{2}$ \quad Center for Research and Development in Mathematics and Applications (CIDMA), \linebreak
Department of Mathematics, University of Aveiro, 3810-193 Aveiro, Portugal}

\corres{Correspondence: delfim@ua.pt; Tel.: +351-234-370-668}

\firstnote{These authors contributed equally to this work.}

\abstract{In this article, we develop a simple mathematical GNU Octave/MATLAB code that is easy to modify 
for the simulation of mathematical models governed by fractional-order 
differential equations, and for the resolution of fractional-order 
optimal control problems through Pontryagin's maximum principle 
(indirect approach to optimal control). For this purpose, a fractional-order 
model for the respiratory syncytial virus (RSV) infection is considered. 
The model is an improvement of one first proposed 
by the authors in [Chaos Solitons Fractals 117 (2018), 142--149]. 
The initial value problem associated with the RSV infection fractional model 
is numerically solved using Garrapa's {\texttt{fde12}} solver and two simple methods 
coded here in Octave/MATLAB: the fractional forward {Euler's} method and the 
predict-evaluate-correct-evaluate (PECE) method of Adams--Bashforth--Moulton. 
A fractional optimal control problem is then formulated having treatment as the control. 
The fractional Pontryagin maximum principle is used to characterize the fractional optimal 
control and the extremals of the problem are determined numerically through the implementation 
of the forward-backward PECE method. The implemented algorithms are available on GitHub and, 
at the end of the paper, in appendixes, both for the uncontrolled initial value problem as well 
as for the fractional optimal control problem, using the free GNU Octave 
computing software and assuring compatibility with~MATLAB.}



\begin{document}

\section{Introduction}

{In 1695, L'H\^{o}pital asked Leibniz in a letter about
the possible meaning of a derivative of order $1/2$~\cite{leibniz1849letter}. 
This episode is considered the kilometer zero of the Fractional Calculus road. 
In recent years, the modeling of real-phenomenon with fractional-order 
derivatives has caught the attention of many researchers.  
The associated problems have been modeled and studied using fractional-order 
derivatives to better understand their dynamics.
For problems that arise in biology, ecology, engineering, physics, 
and some other fields of applied sciences, see, e.g., 
\cite{MR1658022,math10214023,SHAH202211211,SINTUNAVARAT202265}.}

{Mathematical models can predict the evolution of an infectious disease,
show the predictable result of an epidemic, and support public health
possible interventions. Compartmental models serve as a basic 
mathematical structure in epidemiology to comprehend the dynamics of such systems.
In the simplest case, the  compartments divide the population into two health states:
susceptible to the infection of the pathogen agent, usually denoted by $S$, and infected by the
pathogen, usually denoted by the letter $I$. 
Phenomenological assumptions explain  
the way that these compartments interact, and the model is constructed from there.
Usually, these models are investigated through systems of ordinary differential equations.
Other compartments could be included. Depending on the disease, the recovered/immune/removed 
compartment, usually denoted by $R$, is very common. To give further realism to the mathematical 
models and consider the influence of the past on the current and future state of the disease, 
recently, fractional order differential equations have been considered. In this regard,  
one can find, for example, research on dengue, Ebola, tuberculosis, and HIV/AIDS
\cite{area2015mathematical,MR3743014,MR3691319}.}

Respiratory syncytial virus (RSV) is a prominent cause of acute lower respiratory 
infection in young children. {Consequently}, RSV is a considerable burden 
on healthcare~systems. 

In a recent study of RSV in Portugal~\cite{bandeira2022burden}, it is shown that  
RSV is accountable for a substantial number of hospitalizations in children, 
especially when they have less than one year old. Hospitalizations are mainly motivated  
by healthy children. The authors of~\cite{bandeira2022burden} conclude their study 
claiming that the creation of {a} universal RSV surveillance system to guide prevention 
strategies {are} crucial.

In another context, a surveillance system was already implemented in Florida in 1999, 
to support clinical decision-making for the prophylaxis of premature newborns. Since this 
infection is seasonal, a local periodic SEIRS mathematical model was
proposed in~\cite{rosa:delfim2018parameter} to describe real data collected by Florida's system.~Later, a nonlocal fractional (non-integer order) model was proposed in~\cite{MR3872489},
where a fractional optimal control problem was formulated and solved.

In this work, we start by introducing dimension corrections to the SEIRS-$\alpha$ epidemic model 
presented in~\cite{MR3872489}. Afterwards, we apply fractional optimal control to the model
having treatment as the control variable. Differently from previous works, our computer codes 
are presented {in} the text and they are intentionally easy to modify in order to motivate 
readers to use them and adapt them to their own models and to their own contexts. 

When $\alpha=1$, a fractional compartmental model represents a classical compartmental model. 
Therefore, the presented codes can also solve classical optimal control models; 
although, in that case, we suggest the reference~\cite{MR4091761} 
as a preferential option in such a scenario.  

By providing the code of algorithms in an open programming language, we believe that 
our work contributes to reducing the alleged ``replication crisis'' in science
and, in particular, in the field of dynamic optimization and control in biomedical research.
This is a crisis due to the fact that many scientific studies are difficult or 
impossible to validate through replication~\cite{an2018crisis}. 

{The organization of the paper is as follows. In Section~\ref{sec:2},
we introduce the fractional-order RSV model,
correcting the model first presented in~\cite{MR3872489}.
The numerical resolution of the fractional RSV model is presented by
three algorithms, being the subject of Section~\ref{sec:3}. The fractional optimal control 
of RSV transmission is the subject of Section~\ref{sec:4}. We end with 
conclusions and possible future work in Section~\ref{sec:5}.}

\section{A Fractional-Order RSV Model}
\label{sec:2}

We consider that the population under study consists of susceptible ($S$), 
infected but not yet infectious ($E$), infected and infectious ($I$), 
and recovered ($R$) individuals. A characteristic feature of RSV is that immunity 
after infection is temporary, so the recovered individuals become susceptible 
again~\cite{Weber2001}. Let parameter $\mu$ denote the birth rate, which we assume 
equal to the mortality rate; individuals leave the latency period and become infectious 
at a rate $\varepsilon$; $\gamma$ be the rate of loss of immunity; and $\nu$ be 
the rate of loss of infectiousness. We assume the latency period to be
equal to the time between infection and the first symptoms. The influence 
of seasonality on the transmission parameter $\beta$ is modeled by the
\emph{cosine} function. As in~\cite{zhang2012existence}, we consider that the annual
recruitment rate is seasonal due to school opening/closing periods.
Our system of fractional differential equations, the SEIR-$\alpha$ 
epidemic model presented in~\cite{MR3872489} is given by 

\vspace{6pt}
\begin{equation}
\label{eq:modSEIRS0}
\begin{cases}
\Dl S(t) =  \lambda(t)-\mu S(t)-\beta(t) S(t) I(t) +\gamma R(t),\\[1mm]
\Dl E(t) = \beta(t) S(t) I(t) -\mu E(t)-\varepsilon E(t),\\[1mm]
\Dl I(t) =\varepsilon E(t)-\mu I(t)-\nu I(t),\\[1mm]
\Dl R(t) =\nu I(t)-\mu R(t)-\gamma R(t),
\end{cases}
\end{equation}
where  $\beta(t)=b_0(1+b_1\cos(2 \pi t +\Phi))$ and $\Dl$
denote the left Caputo derivative of order $\alpha \in (0,1]$
\cite{MR1658022}. The parameter $b_0$ is the mean
of the transmission parameter and $b_1$ is the amplitude 
of the seasonal fluctuation in the transmission parameter, $\beta$. 
Here, \mbox{$\lambda(t)=\mu(1 + c_1  \cos( 2 \pi t +\Phi) )$} is the recruitment rate
(including newborns and immigrants), where parameter $c_1$ is the amplitude
of the seasonal fluctuation in the recruitment parameter, $\lambda$,
while $\Phi$ is an angle that is chosen in agreement with real data.
Note that in the particular case of $\alpha = 1$ we obtain from 
\eqref{eq:modSEIRS} the SEIRS model of~\cite{UBI:UA}.

Equations of model \eqref{eq:modSEIRS0} do not have appropriate time dimensions. 
Indeed, on the left-hand side the dimension is (time)$^{-\alpha}$, 
while on the right-hand side the dimension is (time)$^{-1}$. We can conclude  
that model \eqref{eq:modSEIRS0} is only consistent when $\alpha=1$. 
For more details about the importance of consistency of dimensions, 
we refer the reader to~\cite{MR3808497,MR3928263} 
and the references therein. Hence, we correct system \eqref{eq:modSEIRS0} as follows:
\begin{equation}
\label{eq:modSEIRS}
\begin{cases}
\Dl S(t) =  \lambda(t)-\mu^\alpha S(t)-\beta(t) S(t) I(t) +\gamma^\alpha R(t),\\[1mm]
\Dl E(t) = \beta(t) S(t) I(t) -\mu^\alpha E(t)-\varepsilon^\alpha E(t),\\[1mm]
\Dl I(t) =\varepsilon^\alpha E(t)-\mu^\alpha I(t)-\nu^\alpha I(t),\\[1mm]
\Dl R(t) =\nu^\alpha I(t)-\mu^\alpha R(t)-\gamma^\alpha R(t),
\end{cases}
\end{equation}
where  $\beta(t)=b_0^\alpha(1+b_1\cos(2 \pi t +\Phi))$ 
and $\lambda(t)=\mu^\alpha(1 + c_1  \cos( 2 \pi t +\Phi))$.

\section{Numerical Resolution of the Fractional RSV Model}
\label{sec:3}

In this section, we consider an initial value problem that consists 
{of} the fractional system \eqref{eq:modSEIRS} and the following 
initial conditions in terms of percentage of total population:
\begin{equation}
\label{initialSol}
{S(0)=0.426282, \enspace E(0)= 0.0109566,  
\enspace I(0)=0.0275076, \enspace  R(0)= 0.535254.}
\end{equation}
The values of \eqref{initialSol} correspond to the endemic equilibrium 
of the fractional system \eqref{eq:modSEIRS}. Note that because we 
have introduced the dimension correction into the initial model, the resulting 
model differs from the one presented in~\cite{MR3872489}. 

The RSV model parameters are presented in Table~\ref{tab:param}.
The parameter values $\varepsilon$, $\nu$, and $\gamma$ were obtained 
from~\cite{Weber2001}. The birth rate, $\mu$, is borrowed from~\cite{flhealthcharts}
for the state of Florida. The birth rate is assumed equal to the mortality rate, 
which results in a constant population during the duration of the study. Analogously 
to~\cite{rosa:delfim2018parameter}, the fractional model was fitted to the data of the 
State of Florida, not including the north region, between September 2011
and July 2014. The data was collected from the Florida Department of Health~\cite{flhealth}.
In that process, values of the following parameters were determined by fitting the model: 
(i) the mean of the transmission parameter, $b_0$; (ii) and its relative seasonal amplitude, $b_1$. 
As previously in~\cite{rosa:delfim2018parameter}, we assume here that the amplitude
of the seasonal fluctuation in the recruitment parameter, $c_1$, is equal to $b_1$.  
The angle $\Phi$ is assumed to be $\pi/2$. This value allows the initial value 
of the transmission parameter to be the average, $\beta(0)=b_0$, and the initial value 
of the recruitment rate to also be the average, $\lambda(0)=\mu$.

The new fractional model maintains a better fit to real data than the classical model 
(the absolute error is equal to 1716.12, which surpasses the value of 1719.12 corresponding 
to the classical model. For more details see~\cite{MR3872489}. The best adjustment to 
real data is achieved with a derivative order slightly different than the one determined before: 
$\alpha=0.995$. This value is the one we consider in what follows.

\begin{table}[H]
\caption{RSV model parameters that, {with the exception of $b_0$, $b_1$ and $c_1$}, 
are borrowed from~\cite{rosa:delfim2018parameter} and references cited therein.}
\label{tab:param}
\begin{tabularx}{\textwidth}{CCCCCCCC}\toprule
\boldmath{$\mu$} & \boldmath{$\nu$} & \boldmath{$\gamma$} & \boldmath{$\varepsilon$} & \boldmath{$b_0$}
& \boldmath{$b_1$ }& \boldmath{$c_1$} & \boldmath{$\Phi$} \\ \midrule
 0.0113 & 36 & 1.8 & 91  & {85} & {0.167} & {0.167} & $\pi/2$ \\
\bottomrule
\end{tabularx}
\end{table}

Algorithms designed to obtain the numerical solution of the initial value problem,   
\eqref{eq:modSEIRS} and \eqref{initialSol}, are now implemented under the free GNU Octave software 
(version 7.3.0), a high-level language primarily intended for numerical computations. 
Octave uses a language that is mostly compatible with MATLAB, being free~\cite{octave}. 
In that regard, two known numerical techniques are implemented: the forward {Euler's}
and the predict-evaluate-correct-evaluate (PECE) methods. The obtained solutions 
are compared with the ones obtained through the known and freely available \texttt{fde12} routine.  

\subsection{The {\emph{fde12}} Solver}

Currently, neither GNU Octave nor MATLAB installation includes a built-in 
routine dedicated to the resolution of nonlinear differential equations 
of fractional order. Nevertheless, in the ``MATLAB Central File Exchange'', 
there exists a routine, named \texttt{fde12},  whose implementation, based 
on Adams--Bashforth--Moulton scheme, is due to Garrapa~\cite{fde12}. 
This Octave/MATLAB routine solves fractional order differential equations 
in the Caputo sense. Convergence and accuracy of the numerical method are available  
in~\cite{diethelm2004detailed}. The stability properties of the algorithm 
implemented by \texttt{fde12} are studied in~\cite{garrappa2010linear}. 

{The initial value problem of Equation \eqref{eq:modSEIRS}, 
with initial conditions \eqref{initialSol}, can be solved with 
the \texttt{fde12} function employing 
the implementation available in Appendix~\ref{ap:Cod:fde12}.}

The solution through Garrapa's routine can then be obtained 
by introducing the following instructions in the 
GNU Octave interface:

\begin{alltt}
	>> N = 400; alpha = 0.995;
	>> [t,y] = model_SEIRS_fde12(N, alpha)
\end{alltt}

\subsection{Fractional Forward Euler's Method}

Let us consider the initial value problem (IVP):
\begin{equation}
\label{sode}
\begin{cases}
\Dl y(t)=f(t,y(t)),\; 0<\alpha\leqslant 1,\\
y(0)=y_0,\;  0< t\leqslant t_f,
\end{cases}
\end{equation}
{where $f (t,y(t))$ is a given function that satisfies 
some smooth conditions~\cite{li2015numerical}.}

The function $y(t)$, known as the exact solution, that satisfies the IVP \eqref{sode} 
is not what we will obtain with this procedure. Instead, an approximation of it 
is computed in as many points as we deem necessary.  

The interval $[0,t_f]$ is subdivided into $n$ subintervals $[t_j,t_{j+1}]$ 
of equal size \linebreak $h=t_f/n$ using the nodes $t_j=jh$, for $j=0,1,\ldots,n$. 

Applying  $\Dlmm$ on both sides of \eqref{sode}, 
we obtain the following equivalent Volterra 
integral equation~\cite{diethelm2004detailed}:
\begin{equation*}
y(t)=y_0+  \Dlmm f(t,y(t)). 
\end{equation*}
According to~\cite{li2015numerical}, $\Dlmm f(t,y(t))$ 
is then approximated by a left fractional rectangular formula in such a way
\begin{equation*}
y(t_{n+1})=y_0+\frac{h^\alpha}{\Gamma(\alpha+1)}\sum_{j=0}^{n}b_{j,n+1}f(t_j,y(j)),   
\end{equation*}
where 
$$
b_{j,n+1}=[(n-j+1)^\alpha-(n-j)^\alpha].
$$
This method was elected here since this approach agrees 
with the fact that fractional derivatives are given by an integral.

Applying the fractional forward Euler's method to approximate 
the four variables of our fractional ordinary differential system 
of equations, we obtain the {GNU Octave routine implementation 
that can be found in Appendix~\ref{ap:Cod:Euler}.}

In Figure~\ref{fig:states:fde12_Euler}, the  solution of the system of fractional 
differential Equation \eqref{eq:modSEIRS}, with initial conditions \eqref{initialSol},
obtained by the \texttt{fde12} solver (solid line) is compared with the correspondent  
solution of the forward {Euler's} method (dashed line). The same discretization is used 
by both implementations in the interval $[0, t_{final}]$, considering 400 knots  
and the step size $h=t_{final}/399$. We can verify that the graphics of system 
variables have strong oscillations. As a consequence, the approximations obtained 
by the forward {Euler's} method, a very simple algorithm, have difficulty following 
the solution of a more sophisticated  and robust implementation like 
the one provided by the \texttt{fde12} solver. 

\vspace{-9pt}
\begin{figure}[H]
\subfloat[\centering{}]{%
\includegraphics[scale=0.405]{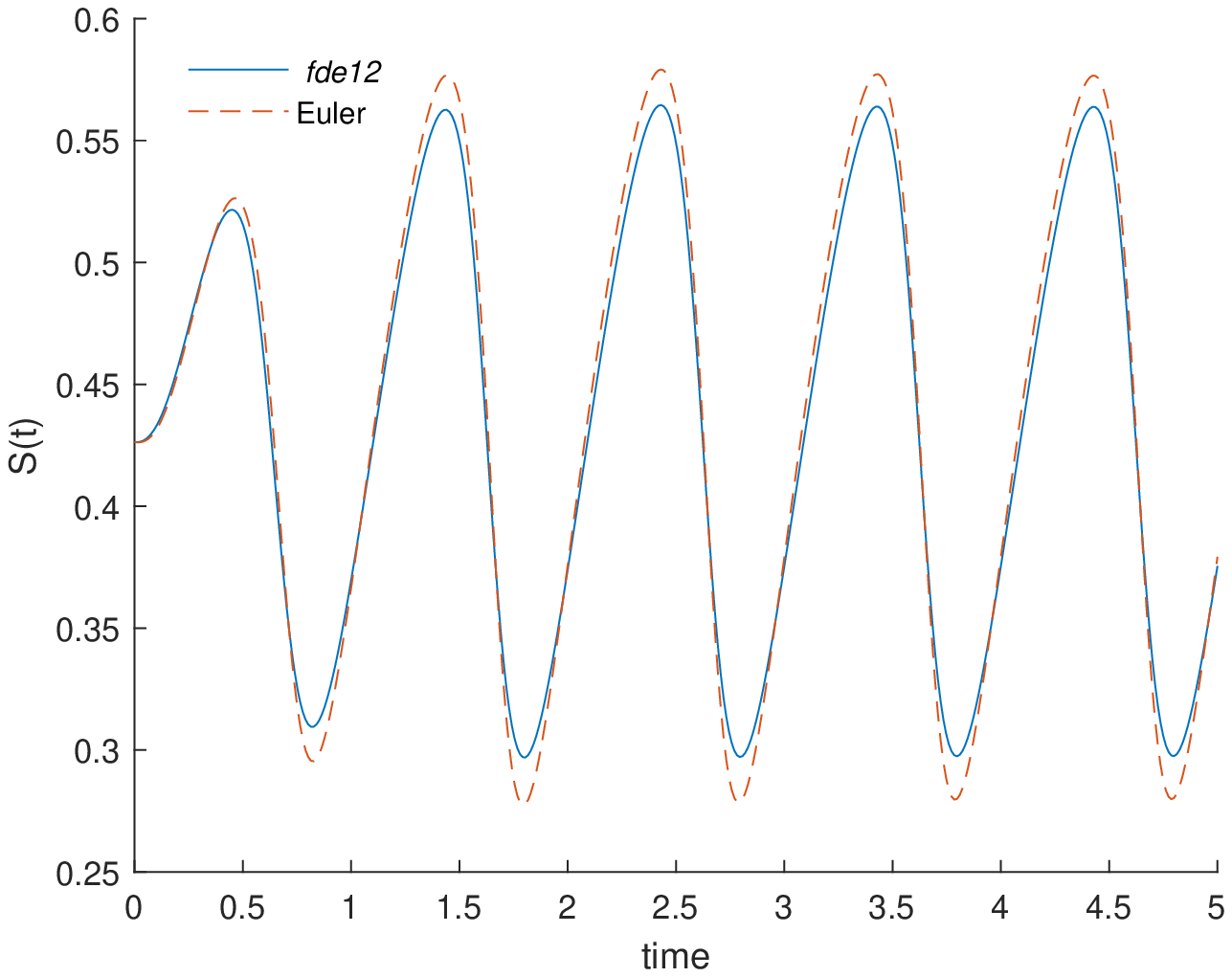}\label{fig:I:seirs_euler}}
\hspace*{0.5cm}
\subfloat[\centering{}]{%
\includegraphics[scale=0.405]{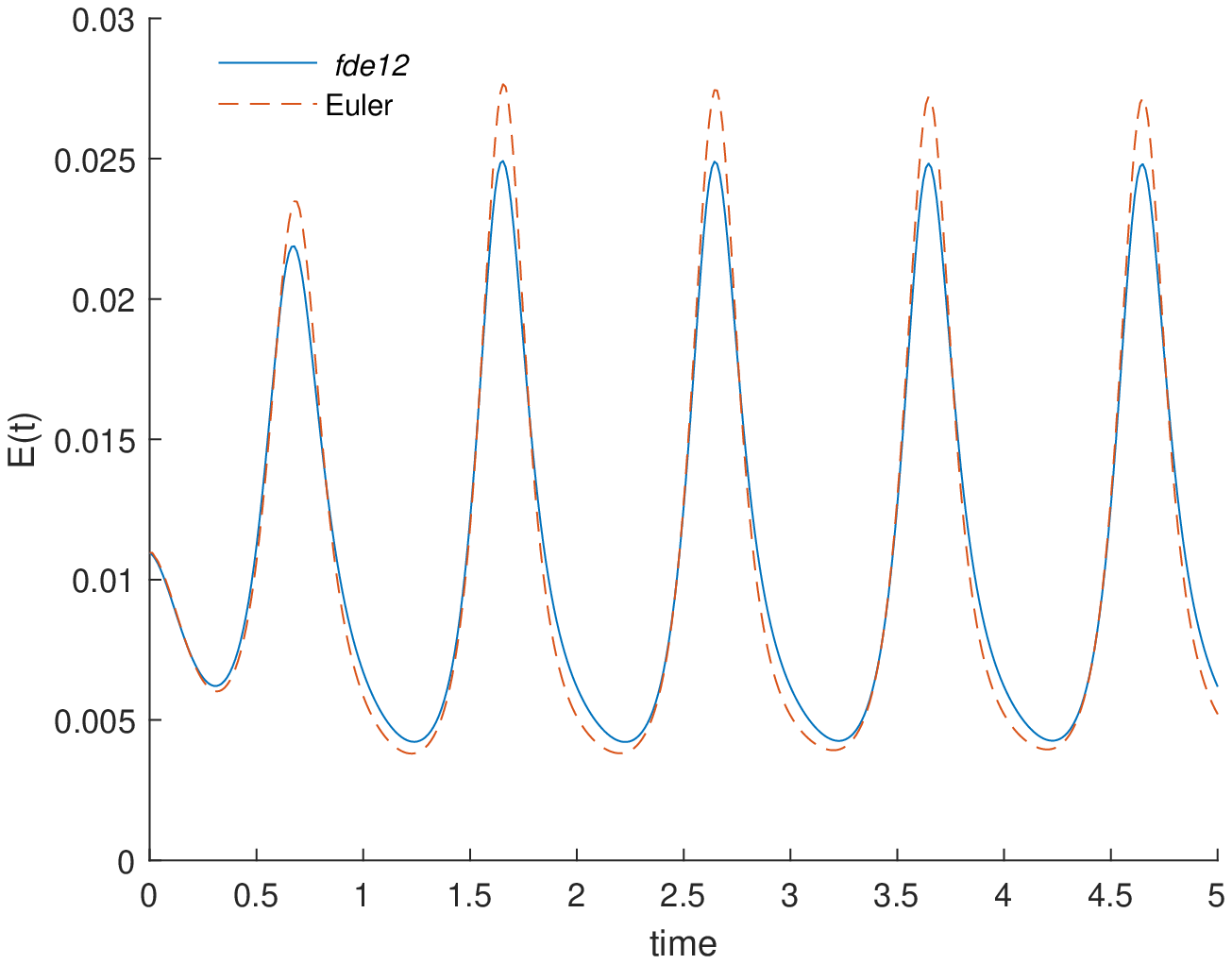}}\\[0.3cm]
\subfloat[\centering{}]{%
\includegraphics[scale=0.405]{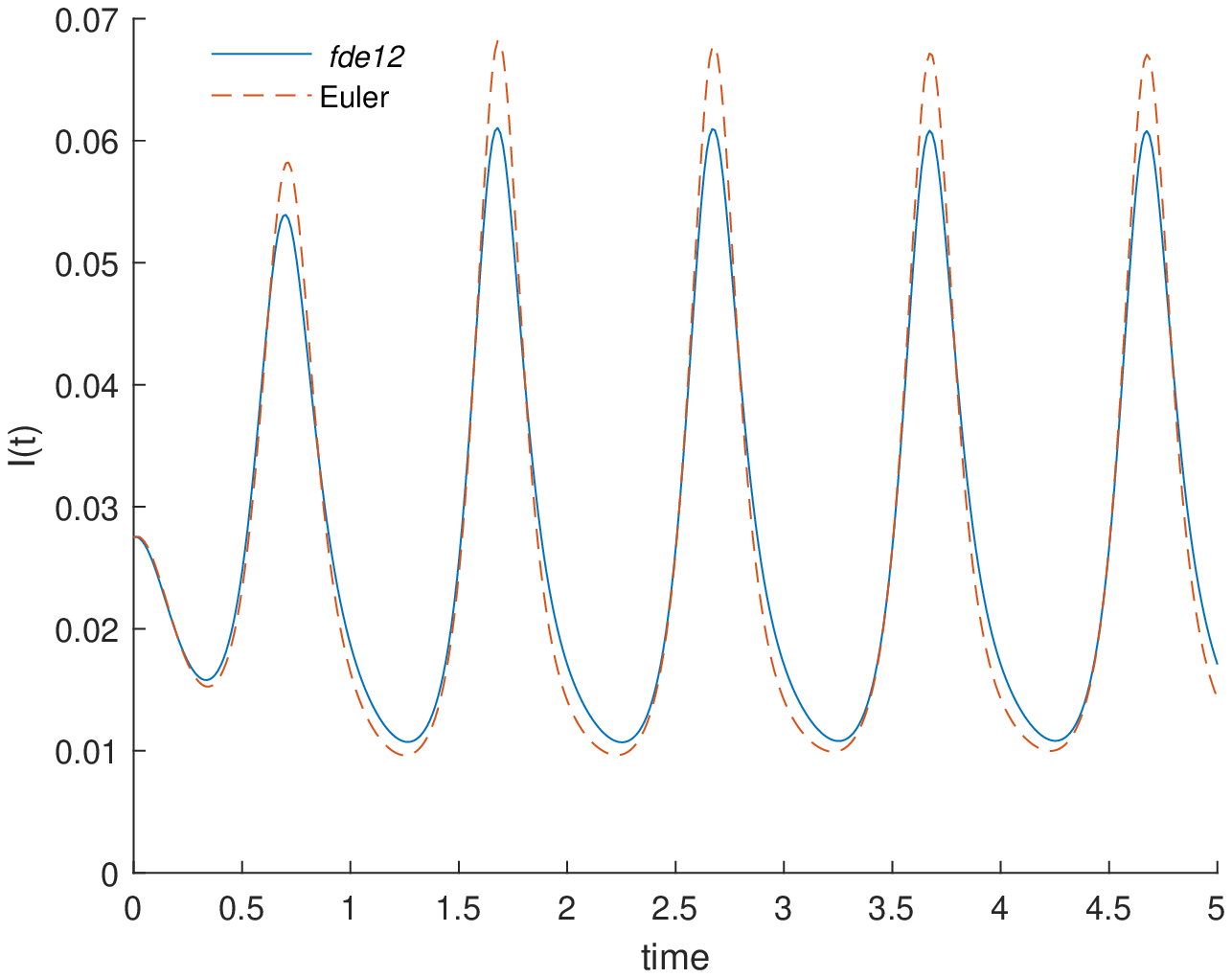}}
\hspace*{0.5cm}
\subfloat[\centering{}]{%
\includegraphics[scale=0.405]{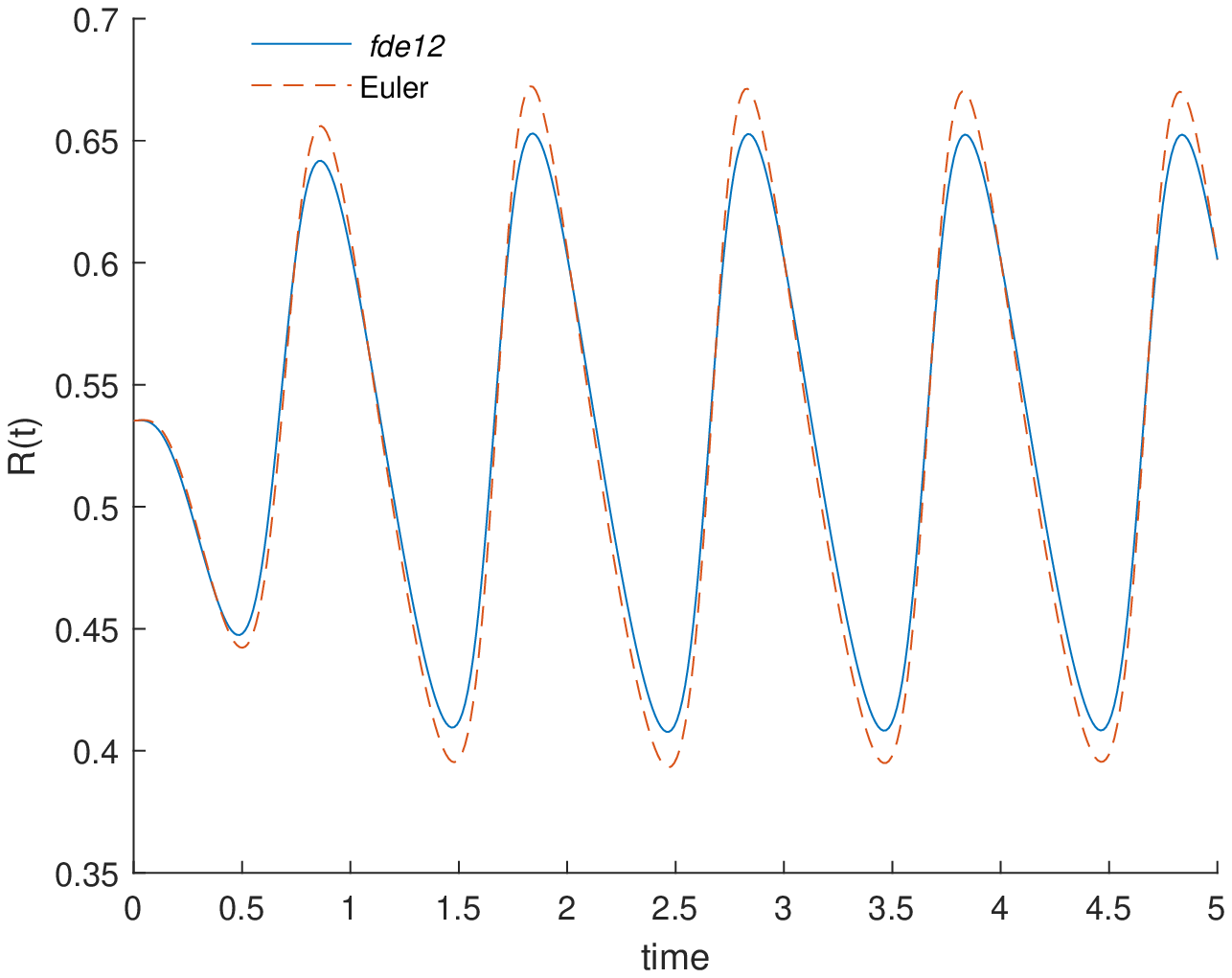}}
\caption{State variables of the fractional differential system \eqref{eq:modSEIRS}, 
considering $\alpha=0.995$, determined with the \texttt{fde12} solver 
and with Euler's method: (\textbf{a}) Variation of the number of susceptible individuals; 
(\textbf{b}) Variation of the number of exposed individuals; 
(\textbf{c}) Variation of the number of infected individuals; 
(\textbf{d}) Variation of the number of recovered individuals.}
\label{fig:states:fde12_Euler}
\end{figure}

The solver \texttt{fde12} is a sophisticated routine that, in some cases, 
considers a number of knots different from the one proposed by the \emph{user}. 
We tested a few numbers of knots and verified  that the number 400 was not changed 
by the solver. This allows the direct comparison with the approximations obtained 
by other routines, dispensing the use of other tools (e.g., interpolation) 
and their associated errors.  

For illustration purposes, we follow the instructions used in the GNU Octave interface  
to obtain Figure~\ref{fig:states:fde12_Euler}a, which exhibits the variation 
of the number of susceptible individuals:

\begin{alltt}
	   figure
	   hold on
	   {plot( t,yf(1,:),t,ye(1,:),'--')}
	   xlabel('time')
	   ylabel('S(t)')
	   legend( '{\textbackslash}it fde12','Euler');
	   legend('boxoff')
	   set(gca,'XTick',[0 0.5 1 1.5 2 2.5 3 3.5 4 4.5 5])
	   hold off
\end{alltt}
where ~{\tt yf(1,:)} is the vector of values of variable $S(t)$ (susceptibles) 
determined by {\tt fde12} solver, while ~{\tt ye(1,:)} is the vector of values 
of variable $S(t)$ determined by Forward Euler's method.

The norm of the difference vector, the absolute difference between the results obtained 
by the \texttt{fde12} solver and the ones obtained by the proposed implementation 
of Euler's method, is presented in Table~\ref{tab:euler} using norms 1, 2, and $\infty$. 
Since Euler's method has a global error of order one~\cite{li2015numerical}, the error
bound depends linearly on the step size, $h$. Therefore, reducing the step size should 
lead to greater accuracy in the approximations.

\begin{table}[H]
\caption{Difference between results of \texttt{fde12} solver and Euler's method 
with norms 1, 2, and $\infty$.}\label{tab:euler}
{\begin{tabularx}{\textwidth}{cCCCC}\toprule
\textbf{Variable} & \boldmath{$S(t)$} & \boldmath{$E(t)$} 
& \boldmath{$I(t)$} &\boldmath{$R(t)$}\\ \midrule
$\|{\it fde12}-\mathrm{Euler}\|_1$ & 3.89849 & 0.297874 & 0.776241 & 3.8815 \\
$\|{\it fde12}-\mathrm{Euler}\|_2$ &  0.221838 & 0.0196 & 0.0512966 & 0.22177   \\
$\|{\it fde12}-\mathrm{Euler}\|_{\infty}$ 
& 0.0197738 & 0.00285878 & 0.00745357 & 0.019802 \\ \bottomrule
\end{tabularx}}
\end{table}

\subsection{PECE Algorithm}
\label{sec:PECE:Alg}

Based on~\cite{diethelm2005algorithms} and the references cited therein,
we now propose an implementation of the predict-evaluate-correct-evaluate (PECE) 
me\-thod of Adams--Basforth--Moulton. The code is relatively simple, easy to modify, 
and can be tailored to solve a particular nonlinear fractional differential model 
with constant, or time-varying, coefficients. 

Applying the PECE method to approximate the four variables 
of our fractional ordinary differential system of equations, 
we obtain the  GNU Octave routine implementation 
that can be found in Appendix~\ref{ap:Cod:PECE}.

In Figure~\ref{fig:states:fde12_pece}, the solution of the system of fractional 
differential Equation \eqref{eq:modSEIRS}, with initial conditions \eqref{initialSol}, 
obtained by the \texttt{fde12} solver (solid line) is compared with 
the corresponding solution of the PECE method (dashed line). 
We observe that the PECE method produces a better approximation than Euler's method since 
both curves in each plot are almost indistinguishable. As before, here the same discretization 
is also used by both implementations in the interval $[0, t_{final}]$, considering 400 knots  
and the step size  $h=t_{final}/399$. 

The norm of the absolute difference between the results obtained by the \texttt{fde12} 
solver and the proposed implementation of the PECE method are presented in Table~\ref{tab:pece} 
making use of norms 1, 2, and $\infty$. Since PECE's method has a global error  
of order two~\cite{Diethelm2002},  the error bound depends quadratically on the step size, $h$. 
This explains why the results obtained from the PECE implementation are better than those
obtained from the forward Euler's~method.

In Section~\ref{sec:4}, a fractional optimal control problem of the model 
is presented. The computation of the corresponding  optimal solution is carried out 
in a forward-backward scheme with the PECE algorithm because \emph{custom made} 
dedicated algorithms are, in general, more efficient and can handle 
more complex models comparatively with generic codes.

\begin{table}[H]
\caption{Difference between results of the \texttt{fde12} solver 
and the PECE implementation with norms 1, 2, and $\infty$.}
\label{tab:pece}
{\begin{tabularx}{\textwidth}{cCCCC}\toprule
\textbf{Variable} & \boldmath{$S(t)$} & \boldmath{$E(t)$} & \boldmath{$I(t)$} &\boldmath{$R(t)$}\\ \midrule
$\|{\it fde12}-\mathrm{PECE}\|_1$ &  0.191041 & 0.0162465 & 0.0415721 & 0.18758\\
$\|{\it fde12}-\mathrm{PECE}\|_2$ & 0.0115686 & 0.00106802 & 0.00269032 & 0.0113171\\
$\|{\it fde12}-\mathrm{PECE}\|_{\infty}$ 
& 0.00133593 & 0.000135793 & 0.000322856 & 0.00128548 \\ \bottomrule
\end{tabularx}}
\end{table}

\begin{figure}[H]
\subfloat[\centering{}]{%
\includegraphics[scale=0.405]{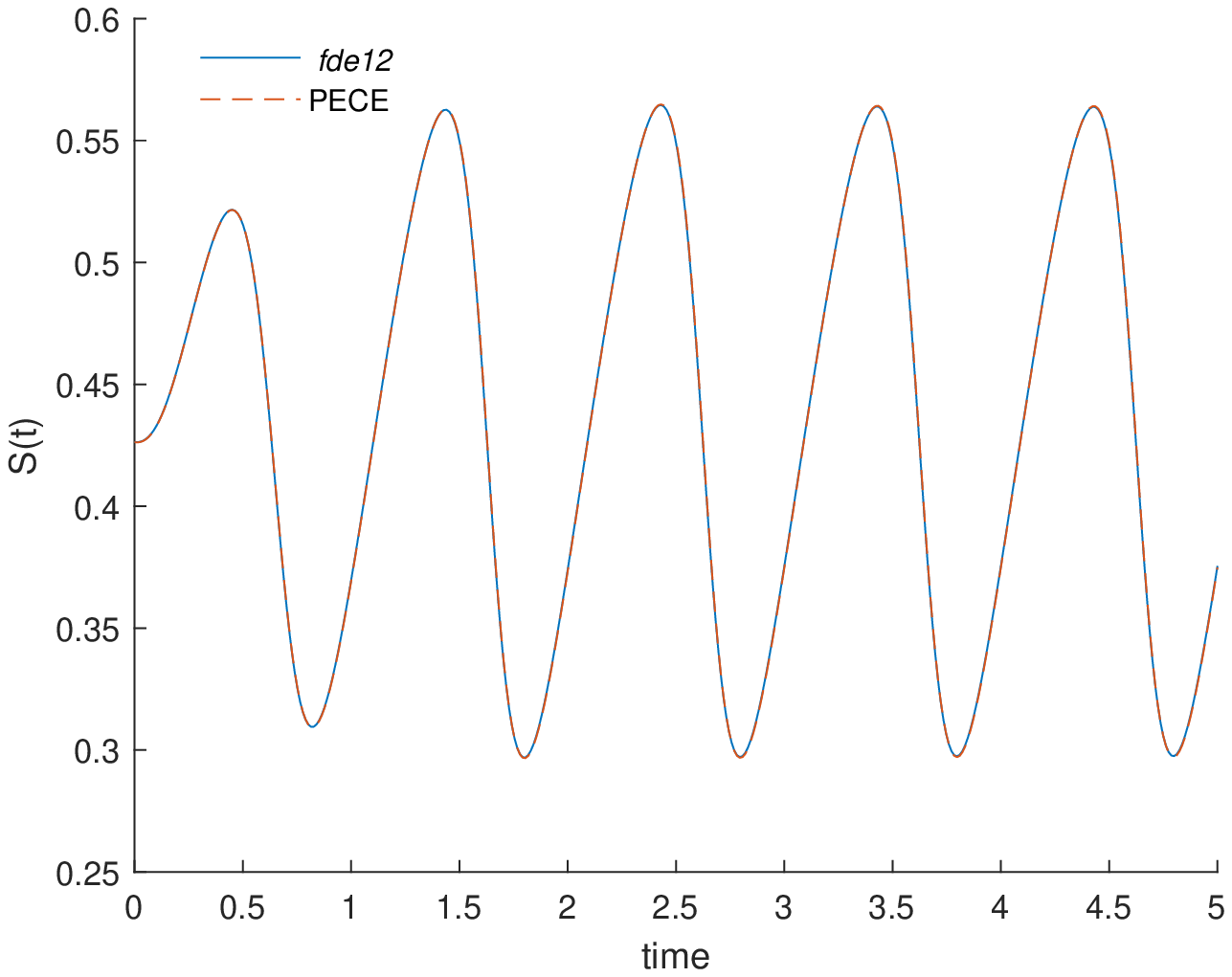}}
\hspace*{0.6cm}
\subfloat[\centering{}]{%
\includegraphics[scale=0.405]{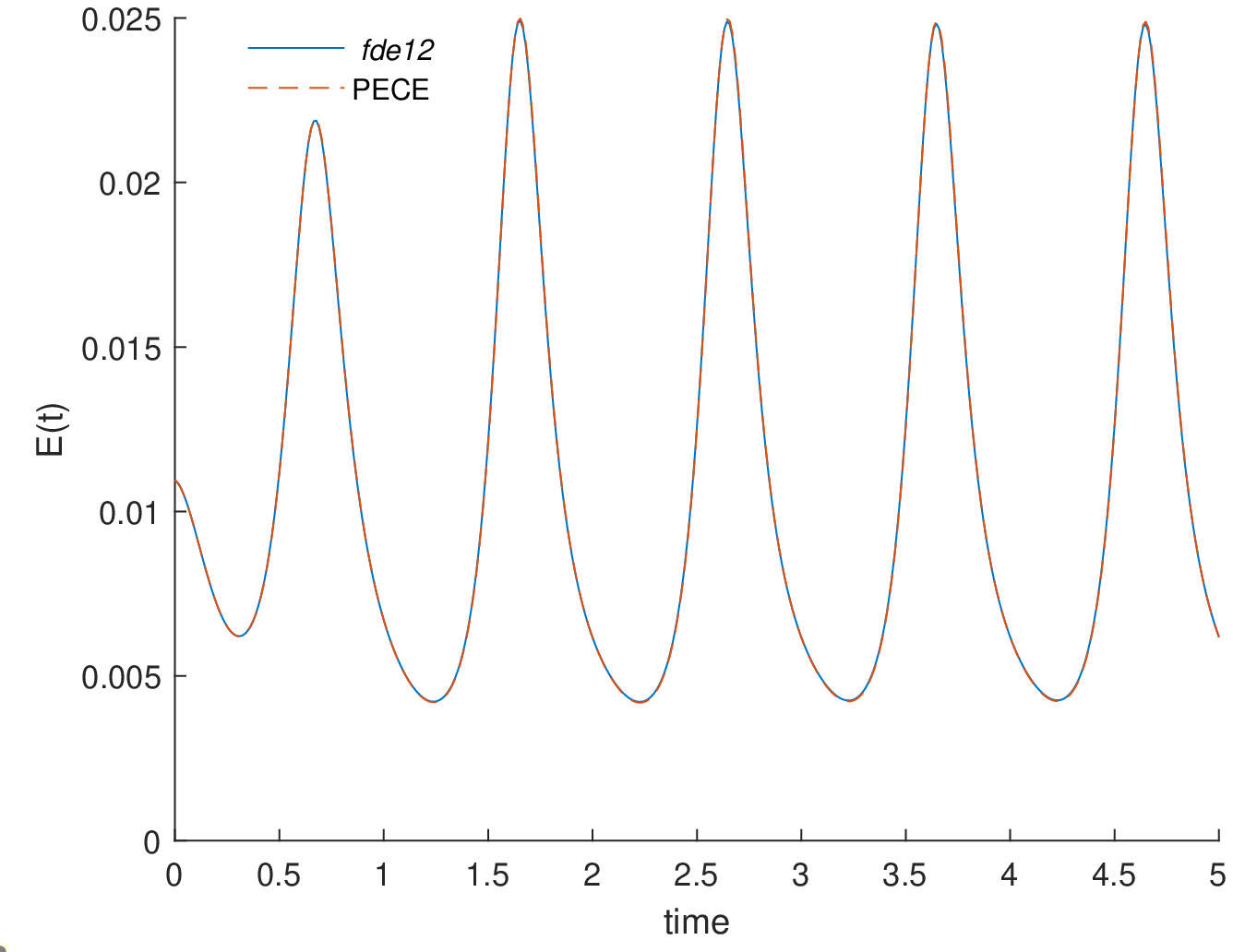}}\\[0.3cm]
\subfloat[\centering{}]{%
\includegraphics[scale=0.405]{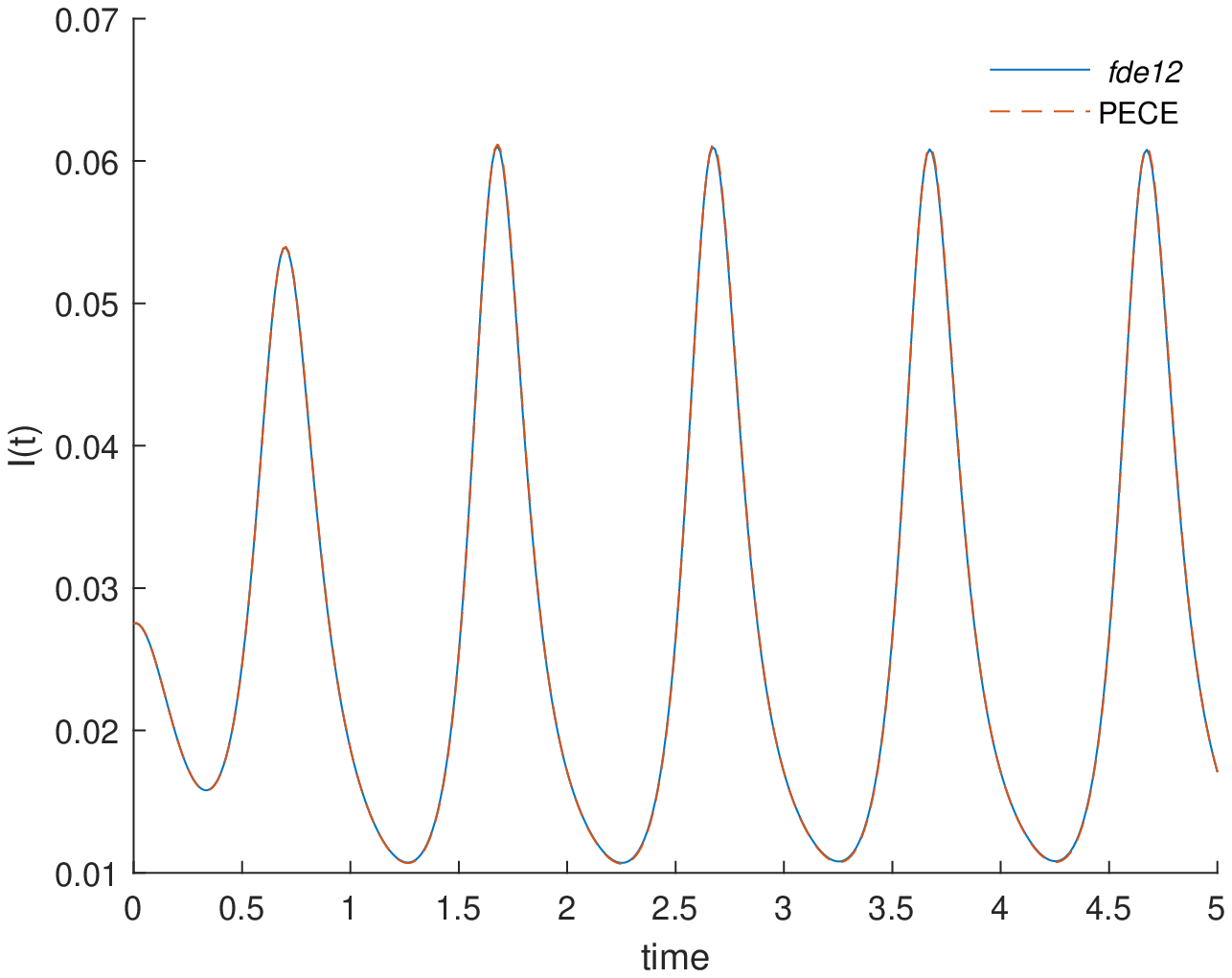}}
\hspace*{0.6cm}
\subfloat[\centering{}]{%
\includegraphics[scale=0.405]{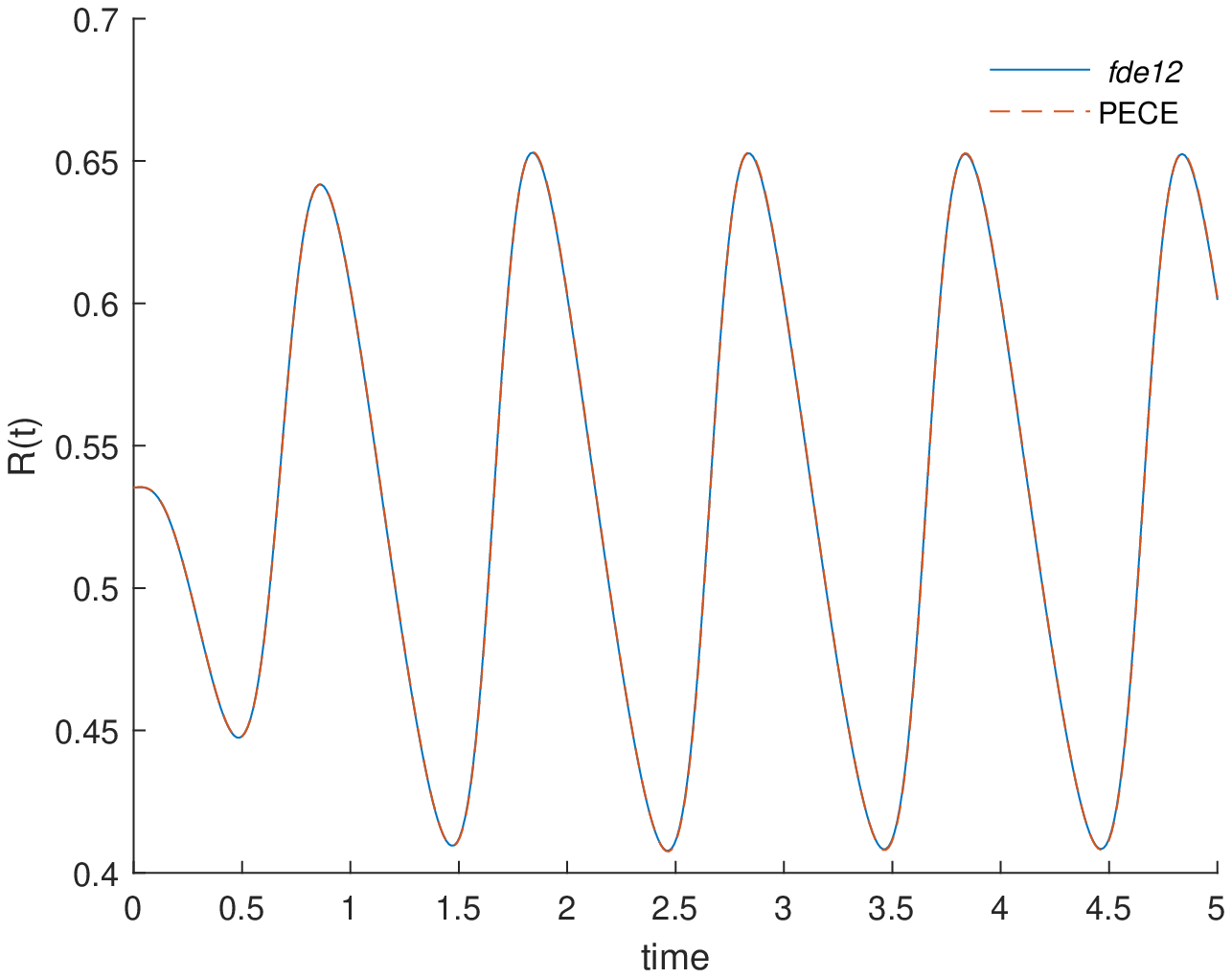}}
\caption{State variables of the fractional differential system \eqref{eq:modSEIRS}, 
considering  {$\alpha=0.995$}, determined with the \texttt{fde12} solver 
and with our PECE implementation: (\textbf{a}) Variation of the number 
of susceptible individuals; (\textbf{b}) Variation of the number of exposed individuals; 
(\textbf{c}) Variation of the number of infected individuals; 
(\textbf{d}) Variation of the number of recovered individuals.}
\label{fig:states:fde12_pece}
\end{figure}

\section{Fractional Optimal Control of RSV Transmission}
\label{sec:4}

The evolution of the variables of the model {depends} on some circumstances 
that can be controlled. In what concerns RSV disease, treatment is commonly 
used as the control due to its relevance in a hospital context (limitation 
in the number of beds and other resources). Hence, we consider the following fractional
optimal control problem: to minimize the number of infectious individuals and the cost
associated to control the disease with the treatment of the patients, that is,
\begin{equation}
\label{cost-functional}
\min ~\mathcal{J}(I,\mathbbm{T})
=\int_0^{t_f} \left(\kappa_1\,I(t)+ \kappa_2\, \mathbbm{T}^2(t)\right) ~dt
\end{equation}
with given $0<\kappa_1,\kappa_2 <\infty$, 
subject to the fractional control system
\begin{equation}
\label{eq:modSEIRS_control}
\begin{cases}
\Dl S(t) =  \lambda(t)-\mu^\alpha S(t)-\beta(t) S(t) I(t) +\gamma^\alpha R(t),\\[1mm]
\Dl E(t) = \beta(t) S(t) I(t) -\mu^\alpha E(t)-\varepsilon^\alpha E(t),\\[1mm]
\Dl I(t) =\varepsilon^\alpha E(t)-\mu^\alpha I(t)-\nu^\alpha I(t)-\mathbbm{T}(t) I(t),\\[1mm]
\Dl R(t) =\nu^\alpha I(t)-\mu^\alpha R(t)-\gamma^\alpha R(t)+\mathbbm{T}(t) I(t)
\end{cases}
\end{equation}
and given initial conditions
\begin{equation}
\label{ocp:ic}
S(0),E(0),I(0),R(0)\geqslant 0.
\end{equation}
Here, $\mathbbm{T}$ is the control variable, which designates \emph{treatment}.
Note that in absence of treatment, that is, for $\mathbbm{T}(t) \equiv 0$,
then the control system \eqref{eq:modSEIRS_control} reduces to the SEIRS-$\alpha$
dynamical system \eqref{eq:modSEIRS}. The set of admissible control functions is
\begin{equation}
\label{Omega:set}
\Omega=\left\{\mathbbm{T}(\cdot)\in L^{\infty}(0,t_f):
0\leqslant \mathbbm{T}\leqslant \mathbbm{T}_{\max},\forall t\in[0,t_f]\right\}.
\end{equation}

Two approaches can be chalked to solve optimal control problems: direct and indirect. 
In direct methods, the resolution of the fractional optimal control problem 
is performed through the application of a variety of discretization and numerical tools
\cite{MR3854267}. Indirect methods are based on  Pontryagin's maximum principle 
and are more robust, although less widespread in biological applications since 
they are not as easy to solve as direct approaches~\cite{MR3443073}.
In what follows, we show how one can take advantage of Octave/MATLAB 
to solve fractional optimal control problems through Pontryagin's
maximum principle, reducing the problem to the solution of a boundary value problem.

Pontryagin's maximum principle (PMP) for fractional optimal control can be used
to solve the problem~\cite{MR3443073,MR4116679}.
The Hamiltonian of our optimal control problem is
\begin{multline*}
\mathcal{H}
= \kappa_1 I +\kappa_2 \mathbbm{T}^2+p_1(\lambda-\mu^\alpha S-\beta S I +\gamma^\alpha R)
+p_2(\beta S I -\mu^\alpha E-\varepsilon^\alpha E)\\
+p_3(\varepsilon^\alpha E-\mu^\alpha I-\nu^\alpha I-\mathbbm{T}I)
+p_4(\nu^\alpha I-\mu^\alpha R-\gamma^\alpha R+\mathbbm{T}I);
\end{multline*}
the optimality condition of PMP ensures that the optimal control is given by
\begin{equation}
\label{eq:ext:cont}
\mathbbm{T}(t)=\min\left\{\max\left\{0,\dfrac{(p_3(t)-p_4(t)) I(t)}{2
\kappa_2}\right\},\mathbbm{T}_{\max}\right\};
\end{equation}
while the adjoint system asserts that the co-state variables
$p_i(t)$, $i = 1,\ldots, 4$, satisfy
\begin{align}
\label{eq:co_states_system}
\begin{cases}
\DrRL \, p_1(t) = p_1(t) \left(\mu^\alpha+\beta(t)I(t)\right)-\beta(t) I(t) p_2(t),\\[1mm]
\DrRL \, p_2(t) =  p_2(t) \left(\mu^\alpha+\varepsilon^\alpha\right)-\varepsilon^\alpha p_3(t), \\[1mm]
\DrRL \, p_3(t)  =  -\kappa_1+\beta(t) p_1(t) S(t)-p_2(t) \beta(t) S(t)\\[1mm]
\qquad\qquad\quad +p_3(t) \left(\mu^\alpha+\nu^\alpha+\mathbbm{T}(t)\right)
-p_4(t)\left(\nu^\alpha+\mathbbm{T}(t)\right),\\[1mm]
\DrRL \, p_4(t) = -\gamma^\alpha p_1(t)+p_4(t)\left(\mu^\alpha+\gamma^\alpha\right),
\end{cases}
\end{align}
which is a fractional system of right Riemann--Liouville derivatives, 
whose operator is represented by $\DrRL$.
In addition, the following transversality conditions hold:
\begin{equation}
\label{eq:transversality}
_{\scriptscriptstyle t}D^{\alpha-1}_{\scriptscriptstyle t_f}
p_i\bigm|_{\scriptscriptstyle t_f}=0
\Leftrightarrow _{\scriptscriptstyle t}\!I^{1-\alpha}_{\scriptscriptstyle t_f}
p_i\bigm|_{\scriptscriptstyle t_f}=p_i(t_f)=0, \quad i=1,\ldots,4,
\end{equation}
where $_{\scriptscriptstyle t}\!I^{1-\alpha}_{\scriptscriptstyle t_f}$
is the right Riemann--Liouville fractional integral of order $1-\alpha$.

\subsection*{Numerical Resolution of the RSV Fractional Optimal Control Problem}

The optimal control problem \eqref{cost-functional}--\eqref{Omega:set}
is numerically solved with the help of {Pontryagin's} maximum principle 
and its optimality conditions, as discussed at the beginning 
of Section~\ref{sec:4}, implementing a forward-backward
predict-evaluate-correct-evaluate (PECE) method of Adams--Basforth--Moulton
(see Section~\ref{sec:PECE:Alg} for the PECE algorithm). The presented forward-backward 
algorithm generalizes the algorithm proposed in reference~\cite{MR2316829}. 

{First,} we solve system \eqref{eq:modSEIRS_control}
by the PECE procedure with initial conditions for the state variables \eqref{ocp:ic}
given in terms of the percentage of the total population, that is, $S(0) + E(0) + I(0) + R(0) = 1$,
and a guess for the control over the time interval $[0,t_f]$,
and obtain the values of the state variables $S$, $E$, $I$ and $R$.

Applying the change of variable
\[
t'=t_f-t
\]
to the system of adjoint Equation \eqref{eq:co_states_system}
and to the transversality conditions \eqref{eq:transversality},
we obtain the following left Riemann--Liouville fractional initial value problem
\eqref{eq:co_states_fr2}--\eqref{eq:trans2}:
\begin{equation}
\label{eq:co_states_fr2}
\begin{cases}
\DlRL \, p_1(t') = -[p_1(t') \left(\mu^\alpha+\beta(t')I(t')\right)-\beta(t') I(t') p_2(t')],\\[1mm]
\DlRL \, p_2(t') =  -[p_2(t') \left(\mu^\alpha+\varepsilon^\alpha\right)-\varepsilon^\alpha p_3(t')],\\[1mm]
\DlRL \, p_3(t')  =  -[-\kappa_1+\beta(t') p_1(t') S(t')-p_2(t') \beta(t') S(t')\\[1mm]
\qquad\qquad\quad +p_3(t') \left(\mu^\alpha+\nu^\alpha+\mathbbm{T}(t')\right)
-p_4(t)\left(\nu^\alpha+\mathbbm{T}(t')\right)],\\[1mm]
\DlRL \, p_4(t') = -[-\gamma^\alpha p_1(t')+p_4(t')\left(\mu^\alpha+\gamma^\alpha\right)],
\end{cases}
\end{equation}
with initial conditions
\begin{equation}
\label{eq:trans2}
p_i(t')\bigm|_{\scriptscriptstyle t'=0}=0, \quad i=1,\ldots,4.
\end{equation}
In turn, conditions \eqref{eq:trans2} imply that 
$$
\DlRL\, p_i(t') =\Dll p_i(t'),\quad i=1,\ldots,4,
$$
which means that the adjoint system \eqref{eq:co_states_fr2} can be treated 
as a Caputo system of fractional differential equations 
(see, e.g., \cite{MR3443073}, {Section~3.3).}

Given the initial conditions \eqref{eq:trans2},
we solve \eqref{eq:co_states_fr2} with the PECE procedure and
obtain the values of the co-state variables $p_i$, $i=1,\ldots,4$.
The control is then updated by a convex combination of the previous control
and the value from \eqref{eq:ext:cont}. This procedure is repeated
iteratively until the values of the controls at the previous iteration
are very close to the ones at the current iteration.

In our numerical computations, we consider that $\mathbbm{T}_{\max}=1$
and the other parameters are fixed according to Table~\ref{tab:param}.
Such values allow the transmission parameter's initial value to be the
average, $\beta(0)=b_0$, and the recruitment rate initial value
to also be the average, $\lambda(0)=\mu$. Our initial conditions,
given by  \eqref{initialSol}, guarantee the existence of a
nontrivial endemic equilibrium for the system \eqref{eq:modSEIRS_control} 
in the absence of control ($\mathbbm{T}(t) \equiv 0$), corresponding 
to the population system \eqref{eq:modSEIRS} prior introduction of treatment.
Because the World Health Organization's goals for most diseases are usually
fixed for five-year periods, we assumed $t_f=5$.

{The algorithm is implemented in GNU Octave divided into four functions, 
the main function named} \verb|FOCP_PECE|. {The implementation of all 
those functions is available in Appendix~\ref{ap:Cod:FOCP}.}

\begin{figure}[H]
\subfloat[\centering{}]{%
\includegraphics[scale=0.405]{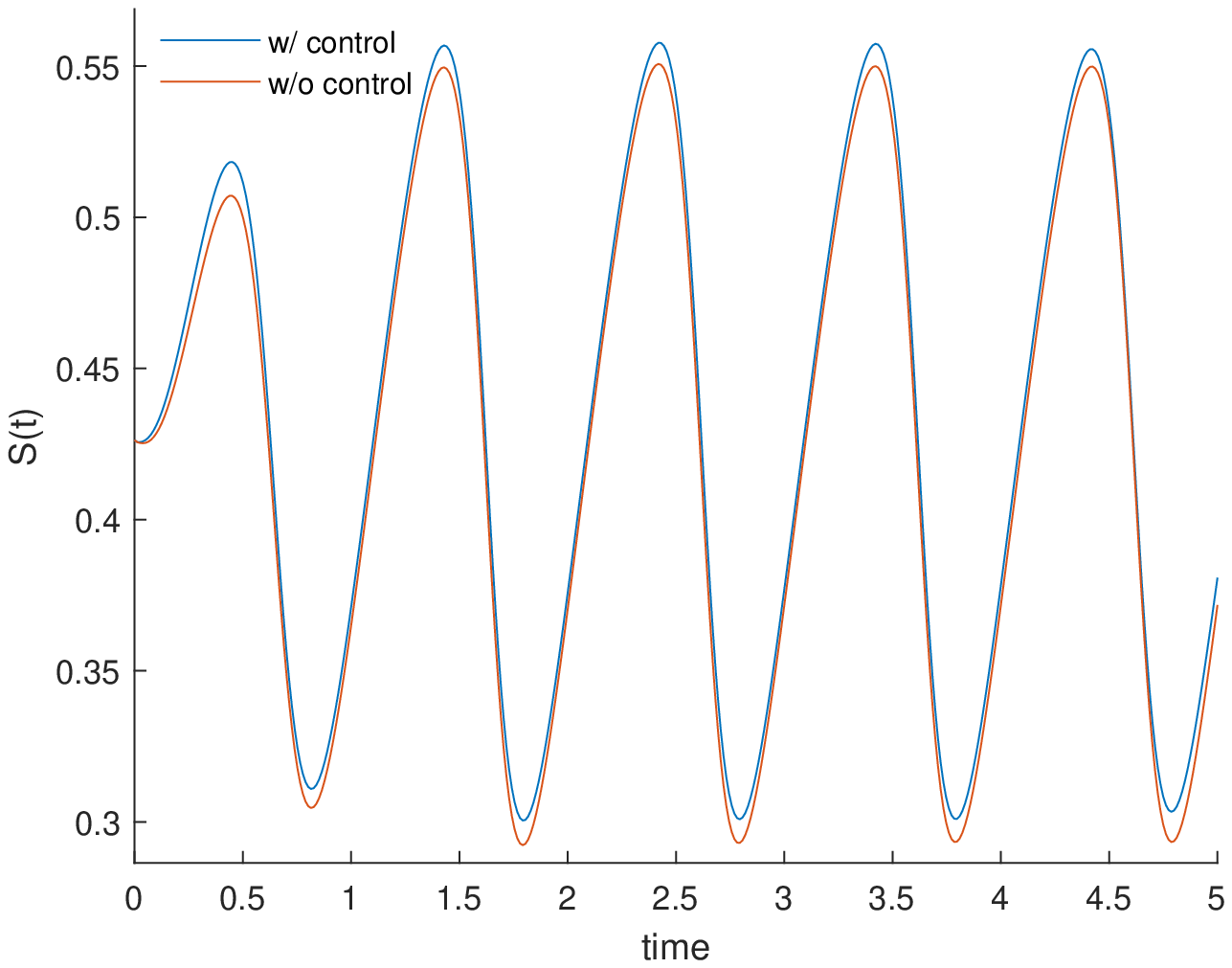}}
\hspace*{0.6cm}
\subfloat[\centering{}]{%
\includegraphics[scale=0.405]{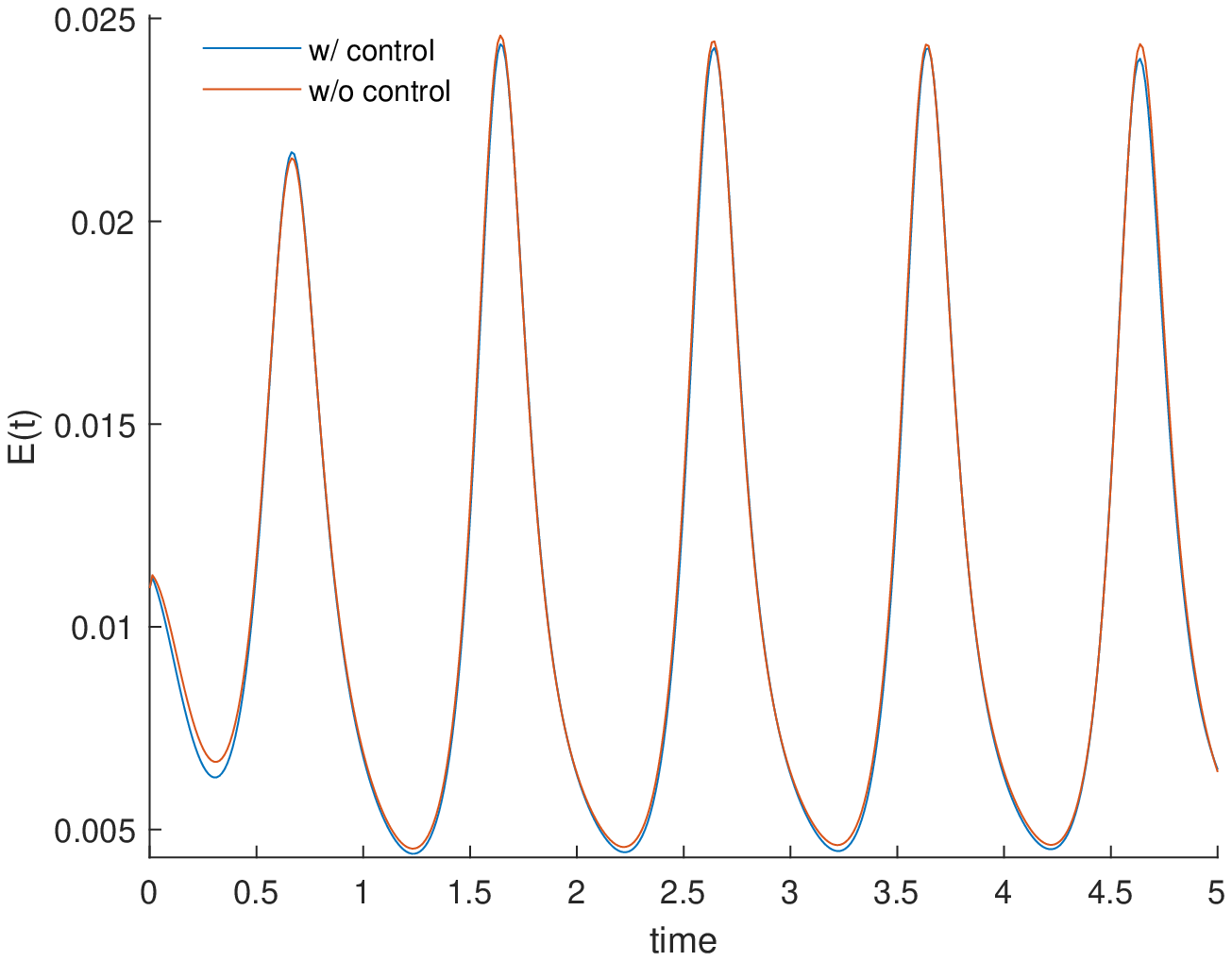}}\\[0.3cm]
\subfloat[\centering{}]{%
\includegraphics[scale=0.405]{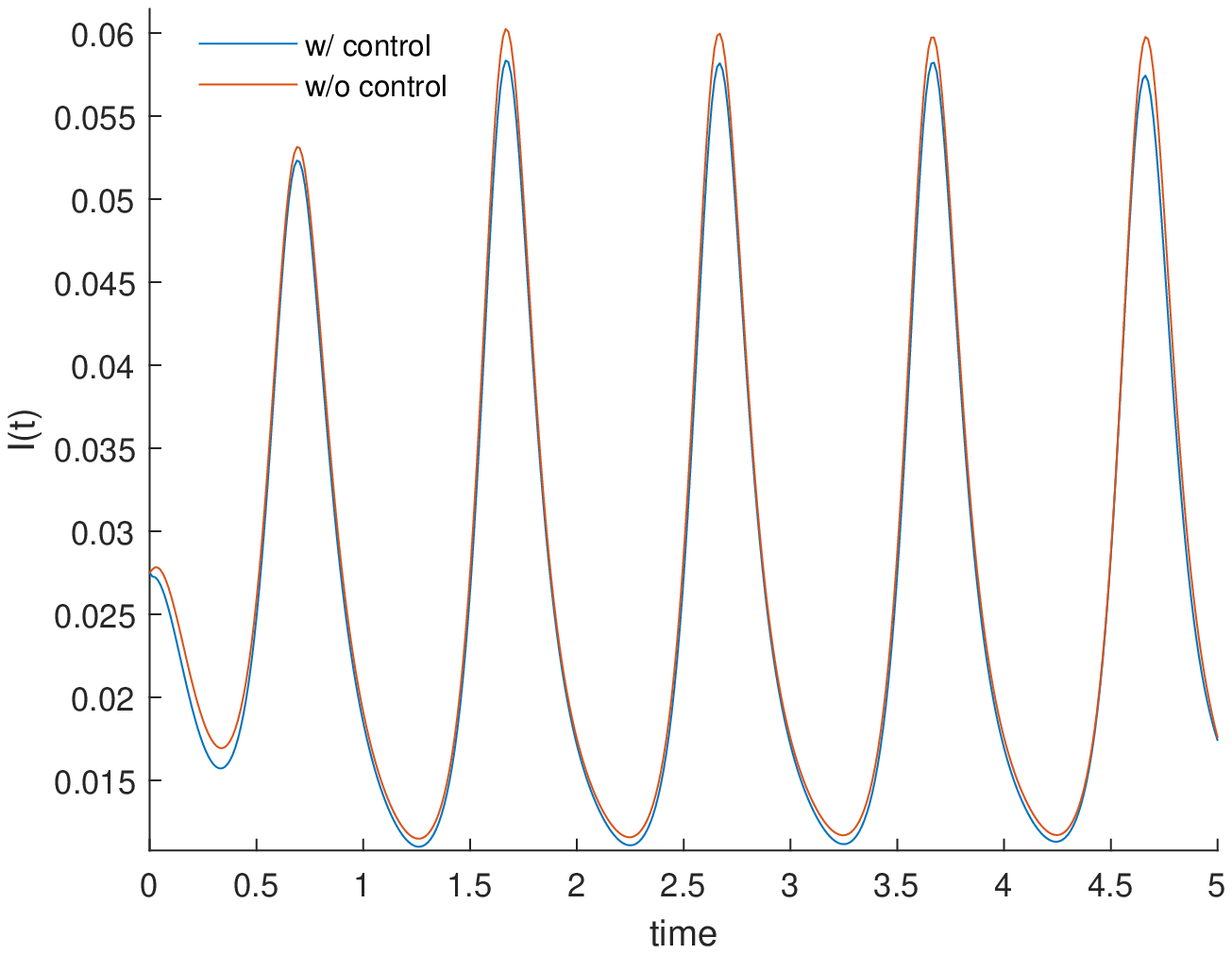}
\label{fig:I:focp_pece}}
\hspace*{0.6cm}
\subfloat[\centering{}]{%
\includegraphics[scale=0.405]{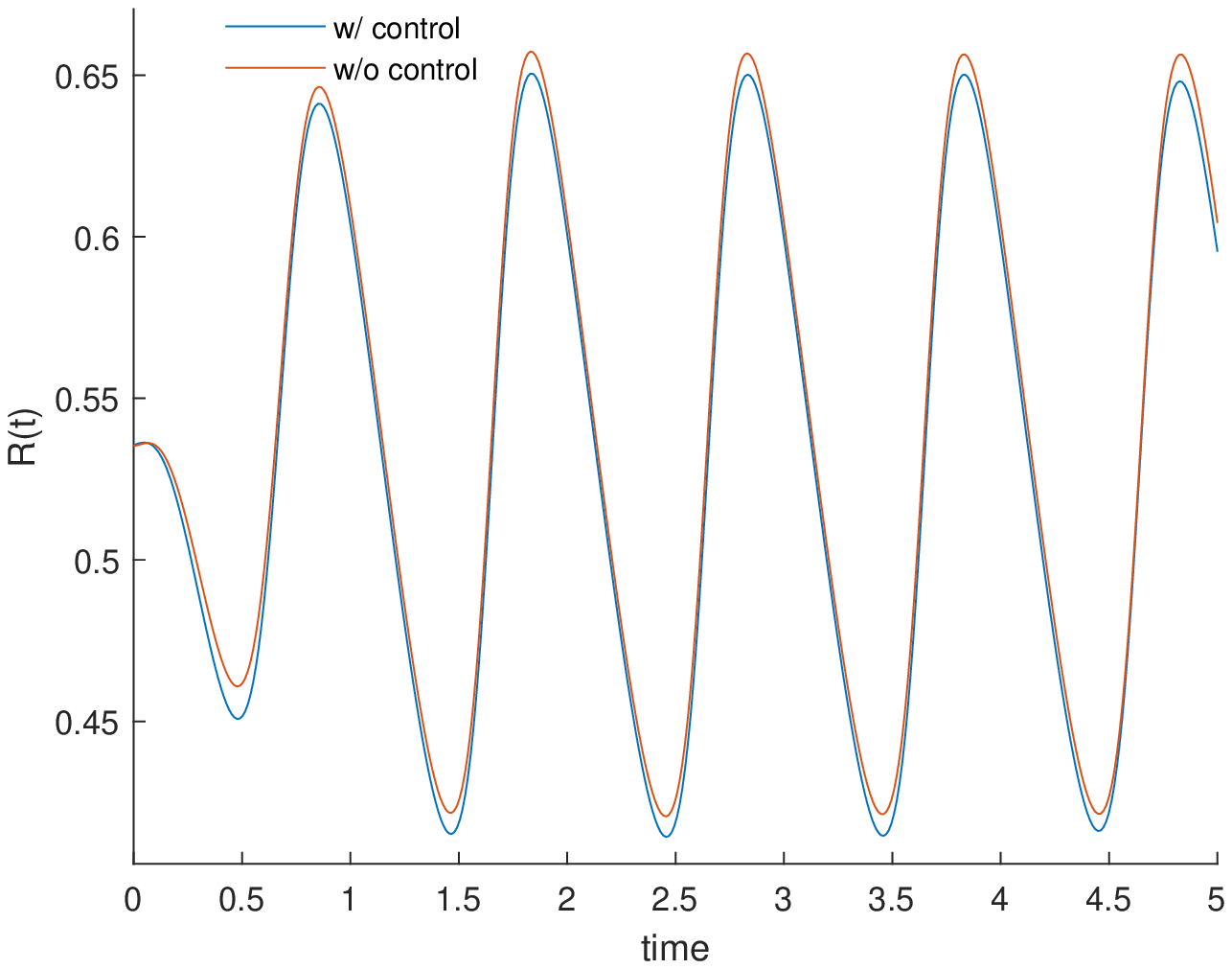}}
\caption{Comparison of optimal state variables for the FOCP, 
defined by \eqref{cost-functional}--\eqref{eq:trans2}, subject 
to the initial conditions \eqref{initialSol} with homonymous 
variables of the original model prior to the use of control treatment. 
(\textbf{a}) Evolution of susceptible individuals. 
(\textbf{b}) Evolution of exposed individuals. 
(\textbf{c}) Evolution of infected individuals. 
(\textbf{d}) Evolution of recovered individuals.}
\label{fig:states:focp_pece}
\end{figure}

The numerical solution of the fractional optimal control problem 
\eqref{cost-functional}--\eqref{eq:trans2} with initial conditions 
\eqref{initialSol}, determined by our \verb|FOCP_PECE| routine 
and associated functions, is exhibited in Figure~\ref{fig:states:focp_pece}. 
The introduction of treatment as control forces a reduction 
in the level of infected individuals as we can see in Figure~\ref{fig:states:focp_pece}c. 

\begin{figure}[H]
\hspace{-12pt}
\includegraphics[scale=0.50]{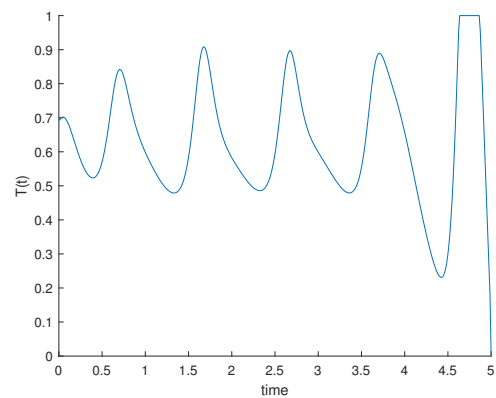}	
\caption{Optimal Control $\mathbbm{T}$ for the RSV fractional 
optimal control problem \eqref{cost-functional}--\eqref{eq:trans2} 
subject to the initial conditions \eqref{initialSol}.}
\label{fig:treatment:focp_pece}
\end{figure}

The evolution of the control treatment is traced in 
Figure~\ref{fig:treatment:focp_pece} and follows 
the seasonality of the RSV infection.

\vspace{-20pt}

\section{Conclusions}
\label{sec:5}

Respiratory syncytial virus (RSV) is the most common cause 
of lower respiratory tract infection in infants and children worldwide.
In addition, RSV causes serious diseases in elderly and immune-compromised 
individuals~\cite{rosa:delfim2018parameter}. In this work, we improved 
a fractional compartmental model for RSV and applied 
optimal control to the resulting model. The Octave/MATLAB codes, developed 
in the computation of numerical solutions, are {available in appendixes} 
and were purposely simplified in order to be easily adapted to other contexts 
and models. We trust this will {motivate more researchers to use fractional 
optimal control in the modeling of} real applications. As future work, 
we plan to investigate the use of others measures to control 
the transmission of the disease and the benefit of the fractional
approach in other contexts and geographical regions.

\vspace{6pt} 

\authorcontributions{Conceptualization, S.R. and D.F.M.T.; 
methodology, S.R. and D.F.M.T.; 
software, S.R. and D.F.M.T.; 
validation, S.R. and D.F.M.T.; 
formal analysis, S.R. and D.F.M.T.; 
investigation, S.R. and D.F.M.T.; 
writing---original draft preparation, S.R. and D.F.M.T.; 
writing---review and editing, S.R. and D.F.M.T.; 
visualization, S.R. and D.F.M.T.
All authors have read and agreed to the published version of the manuscript.} 

\funding{This research was funded by The Portuguese Foundation for Science 
and Technology (FCT -- Funda\c{c}\~ao para a Ci\^encia e a Tecnologia), 
grants number UIDB/50008/2020 (S.R.) and UIDB/04106/2020 (D.F.M.T.).}

\institutionalreview{Not applicable.}

\informedconsent{Not applicable.}
	
\dataavailability{No new data were created.} 

\acknowledgments{{The authors are very grateful 
to three anonymous Reviewers for several important 
comments and suggestions.}}

\conflictsofinterest{The authors declare no conflicts of interest.}

\appendix

\begin{appendices}
\appendixtitles{yes} 
\section{Resolution of the IVP with the {\tt fde12} Function} 
\label{ap:Cod:fde12}

Here, the initial value problem \eqref{eq:modSEIRS} and \eqref{initialSol}
is solved in Octave/MATLAB with the help of the \texttt{fde12} function~\cite{fde12}.

\begin{alltt}
	
	function [t, y] = model_SEIRS_fde12(N,alpha)
		
		% initial conditions
		y0=[0.426282; 0.0109566; 0.0275076; 0.535254];
		
		% Values of parameters
		miu = 0.0113; niu = 36; epsilon = 91; b0 = 85; b1 = 0.167; c1 = 0.167;
		gama = 1.8; tfinal = 5; phi = pi/2;
		ft = linspace(0,tfinal,N); h = tfinal/(N-1); 
		
		% Correction of values of parameters
		miu_ = miu^alpha; niu_ = niu^alpha; epsilon_ = epsilon^alpha;
		gama_ = gama^alpha;
		
		% time-dependent parameters
		flambda = @(t) miu.^alpha.*(1 + c1.* cos( 2.* pi.* t + phi) );
		fbeta = @(t) b0.^alpha.* (1 + b1.* cos( 2.* pi.* t + phi ) );
		
		
		% Differential system of equations of the model 
		fdefun = @(t,y,ft)[flambda(t)-miu_*y(1)-fbeta(t)*y(1)*y(3)+gama_*y(4); ...
		    fbeta(t)*y(1)*y(3)-(miu_+epsilon_)*y(2);
		    epsilon_*y(2)-(miu_+niu_)*y(3); ...
		    niu_*y(3)-miu_*y(4)-gama_*y(4)]; 
		
		% resolution of system with solver fde12 
		[t,y] = fde12(alpha,fdefun,0,tfinal,y0,h,ft);
		
	end
\end{alltt}


{\section{Resolution of the IVP with the Forward Euler Method} 
\label{ap:Cod:Euler}}

Here, the initial value problem \eqref{eq:modSEIRS} and \eqref{initialSol}
is solved in Octave/MATLAB by using the fractional forward Euler's method 
to approximate the four variables of the fractional  
system of equations.

\begin{alltt}
	
	function [t,y] = model_SEIRS_EULER(N,alpha)
		
		% Values of parameters
		miu = 0.0113; niu = 36; epsilon = 91;  gama = 1.8; tfinal = 5; b0 = 85;
		b1 = 0.167; c1 = 0.167; phi = pi/2;
		
		% initial conditions
		S0 = 0.426282; E0 = 0.0109566; I0 = 0.0275076; R0 =  0.535254;
		
		% Correction of values of parameters
		miu_ = miu^alpha; niu_ = niu^alpha; epsilon_ = epsilon^alpha;
		gama_ = gama^alpha;
		
		% time-dependent parameters
		flambda = @(t) miu_*(1 + c1 * cos( 2 * pi * t + phi) );
		fbeta = @(t) b0^alpha.* (1 + b1 * cos( 2 * pi * t + phi ) );
		
		% Initialization of variables
		t = linspace(0,tfinal,N); h = tfinal/N; init = zeros(1,N);
		S = init; E = init; I = init; R = init;
		S(1) = S0; E(1) = E0; I(1) = I0; R(1) = R0;
		beta = fbeta(t); lambda = flambda(t);
		
		
		for j = 2:N
		    aux_s = 0; aux_e = 0; aux_i = 0; aux_r = 0;
		    for k = 1:j-1
		        bk = (j-k+1)^alpha-(j-k)^alpha;
		
		        % Differential system of equations of the model
		        aux_s = aux_s+bk*(lambda(k)-miu_*S(k)-beta(k)*S(k)*I(k) ...
		                +gama_*R(k));
		        aux_e = aux_e+bk*(beta(k)*S(k)*I(k)-(miu_+epsilon_)*E(k));
		        aux_i = aux_i+bk*(epsilon_*E(k)-(miu_+niu_)*I(k));
		        aux_r = aux_r+bk*(niu_*I(k)-miu_*R(k)-gama_*R(k));
		    end
		    S(j) = S0+h^alpha/gamma(1+alpha)*aux_s;
		    E(j) = E0+h^alpha/gamma(1+alpha)*aux_e;
		    I(j) = I0+h^alpha/gamma(1+alpha)*aux_i;
		    R(j) = R0+h^alpha/gamma(1+alpha)*aux_r;
		end
		
		y(1,:) = S; y(2,:) = E; y(3,:) = I; y(4,:) = R;
		
	end
	
\end{alltt}

{\section{Resolution of the IVP with the PECE Method} 
\label{ap:Cod:PECE}}

Now, the initial value problem \eqref{eq:modSEIRS} and \eqref{initialSol}
is solved in Octave/MATLAB by using the predict-evaluate-correct-evaluate
(PECE) method of Adams--Bashforth--Moulton.

\begin{alltt}
	
	function [t,y] = model_SEIRS_PECE(N,alpha)
		
		% Values of parameters
		miu = 0.0113; niu = 36; epsilon = 91;  gama = 1.8; tfinal = 5; 
		b0 = 85; b1 = 0.167; c1=0.167; phi = pi/2;
		
		% initial conditions
		S0 = 0.426282; E0 = 0.0109566; I0 = 0.0275076; R0 =  0.535254;
		
		% Correction of values of parameters
		miu_ = miu^alpha; niu_ = niu^alpha; epsilon_ = epsilon^alpha;
		gama_ = gama^alpha;
		
		% time-dependent parameters
		flambda = @(t) miu_*(1 + c1 * cos( 2 * pi * t + phi) );
		fbeta = @(t) b0^alpha.* (1 + b1 * cos( 2 * pi * t + phi ) );
		
		% Initialization of variables
		t = linspace(0,tfinal,N); h = tfinal/N; init = zeros(1,N);
		beta = fbeta(t); lambda = flambda(t);
		S = init; E = init; I = init; R = init; b = init; a = init;
		S(1) = S0; E(1) = E0; I(1) = I0; R(1) = R0;
		Sp = S; Ep = E; Ip = I; Rp = R;
		
		% computation of coefficients a_k and b_k
		for k = 1:N
		    b(k) = k^alpha-(k-1)^alpha;
		    a(k) = (k+1)^(alpha+1)-2*k^(alpha+1)+(k-1)^(alpha+1);
		end
		
		for j = 2:N
		
		    % First part: prediction
		
		    aux_s = 0; aux_e = 0; aux_i = 0; aux_r = 0;
		    for k = 1:j
		
		        % Differential system of equations of the model
		        aux_s = aux_s+b(j-k+1)*(lambda(k)-miu_*S(k)...
		                -beta(k)*S(k)*I(k)+gama_*R(k));
		        aux_e = aux_e+b(j-k+1)*(beta(k)*S(k)*I(k)...
		                -(miu_+epsilon_)*E(k));
		        aux_i = aux_i+b(j-k+1)*(epsilon_*E(k)-(miu_+niu_)*I(k));
		        aux_r = aux_r+b(j-k+1)*(niu_*I(k)-miu_*R(k)-gama_*R(k));
		    end
		
		    Sp(j) = S0+h^alpha/gamma(1+alpha)*aux_s;
		    Ep(j) = E0+h^alpha/gamma(1+alpha)*aux_e;
		    Ip(j) = I0+h^alpha/gamma(1+alpha)*aux_i;
		    Rp(j) = R0+h^alpha/gamma(1+alpha)*aux_r;
		
		    % Second part: correction
		
		    aux_ss = lambda(j)-miu_*Sp(j)-beta(j)*Sp(j)*Ip(j)+gama_*Rp(j);
		    aux_ee = beta(j)*Sp(j)*Ip(j)-(miu_+epsilon_)*Ep(j);
		    aux_ii = epsilon_*Ep(j)-(miu_+niu_)*Ip(j);
		    aux_rr = niu_*Ip(j)-miu_*Rp(j)-gama_*Rp(j);
		
		    auxx = ((j-1)^(alpha+1)-(j-1-alpha)*j^alpha);
		    aux_s0 = auxx*(lambda(1)-miu_*S(1)-beta(1)*S(1)*I(1)+gama_*R(1));
		    aux_e0 = auxx* (beta(1)*S(1)*I(1)-(miu_+epsilon_)*E(1));
		    aux_i0 = auxx*(epsilon_*E(1)-(miu_+niu_)*I(1));
		    aux_r0 = auxx*(niu_*I(1)-miu_*R(1)-gama_*R(1));
		
		    aux_s = 0; aux_e = 0; aux_i = 0; aux_r = 0;
		    for k = 1:j-1
		
		        % Differential system of equations of the model
		        aux_s = aux_s+a(j-k)*(lambda(k)-miu_*S(k)-beta(k)*S(k)*I(k)...
		                +gama_*R(k));
		        aux_e = aux_e+a(j-k)*(beta(k)*S(k)*I(k)-(miu_+epsilon_)*E(k));
		        aux_i = aux_i+a(j-k)*(epsilon_*E(k)-(miu_+niu_)*I(k));
		        aux_r = aux_r+a(j-k)*(niu_*I(k)-miu_*R(k)-gama_*R(k));
		    end
		
		    S(j) = S0+h^alpha/gamma(2+alpha)*(aux_ss+aux_s0+aux_s);
		    E(j) = E0+h^alpha/gamma(2+alpha)*(aux_ee+aux_e0+aux_e);
		    I(j) = I0+h^alpha/gamma(2+alpha)*(aux_ii+aux_i0+aux_i);
		    R(j) = R0+h^alpha/gamma(2+alpha)*(aux_rr+aux_r0+aux_r);
		end
		
		y(1,:) = S; y(2,:) = E; y(3,:) = I; y(4,:) = R;
		
	end
	
	
\end{alltt}

\vspace{-30pt}

{\section{Numerical Resolution of the Fractional Optimal Control Problem} 
\label{ap:Cod:FOCP}}

Here, we provide our Octave/MATLAB code for the numerical solution 
of the fractional optimal control problem \eqref{cost-functional}--\eqref{eq:trans2} 
with initial conditions \eqref{initialSol}.
 
\begin{alltt}	
	
	
function  [t,y] = FOCP_PECE(N,alpha);
	
	% values assumed as global
	global tfinal miu niu epsilon gama  b0 b1 c1 phi k1 k2 S0 E0 I0 R0;
	
	
	% Values of parameters
	miu = 0.0113; niu = 36; epsilon = 91;  gama = 1.8; tfinal = 5;
	b0 = 85; b1 = 0.167; phi = pi/2; c1 = .167;
	
	% parameters of the algorithm
	k1 = 1; k2 = 0.001; trmax = 1.0; tol = 0.001; test = 1;
	
	
	% initial conditions
	S0 = 0.426282; E0 = 0.0109566; I0 = 0.0275076; R0 =  0.535254;
	
	% initialization of variables
	t = linspace(0,tfinal,N);
	init = zeros(1,N); S = init; E = init; I = init; R = init;
	p1 = init; p2 = init; p3 = init; p4 = init; Ta = init;
	
	% iterations of the numerical method
	
	while test>tol,
	
	    oldS = S; oldE = E; oldI = I; oldR = R;
	    oldp1 = p1; oldp2 = p2; oldp3 = p3; oldp4 = p4; oldTa = Ta;
	
	    % forward PECE iterations
	    [y1] = system1_control(Ta,t,N,alpha);
	    S = y1(1,:); E = y1(2,:); I = y1(3,:); R = y1(4,:);
	
	    % backward PECE iterations
	    [y2] = system2_adjoint(S,I,Ta,t,N,alpha);
	    p1 = y2(1,:); p2 = y2(2,:); p3 = y2(3,:); p4 = y2(4,:);
	
	    % new control
	    Ta = projection((p3-p4).*I/(2*k2),trmax);
	    Ta = ( Ta + oldTa ) / 2;
	
	    % Relative error values for convergence
	    vector = [max(abs(S-oldS))/(max(abs(S))),...
	        max(abs(oldE-E))/(max(abs(E))),...
	        max(abs(oldI-I))/(max(abs(I))),...
	        max(abs(oldR-R))/(max(abs(R))),...
	        max(abs(oldp1-p1))/(max(abs(p1))),...
	        max(abs(oldp2-p2))/(max(abs(p2))),...
	        max(abs(oldp3-p3))/(max(abs(p3))),...
	        max(abs(oldp4-p4))/(max(abs(p4))), ...
	        max(abs(oldTa-Ta))/(max(abs(Ta)))]*100;
	
	    test = max(vector);
	end
	
	y(1,:) = S; y(2,:) = E; y(3,:) = I; y(4,:) = R; y(5,:) = Ta;
	y(6,:) = p1; y(7,:) = p2; y(8,:) = p3; y(9,:) = p4;
	
end
	
	
% function II: resolution of the fractional control system
	
	
function [y]= system1_control(Ta,t,N,alpha)
	
	global  b0 b1 c1 phi miu gama epsilon niu tfinal S0 E0 I0 R0;
	
	% time-dependent parameters
	flambda = @(t) miu^alpha*(1 + c1 * cos( 2 * pi * t + phi) );
	fbeta = @(t) b0^alpha.* (1 + b1 * cos( 2 * pi * t + phi ) );
	
	% Correction of values of parameters
	miu_ = miu^alpha; niu_ = niu^alpha; epsilon_ = epsilon^alpha;
	gama_ = gama^alpha;
	
	% initialization of variables
	beta = fbeta(t); lambda = flambda(t);
	h = tfinal/N; init = zeros(1,N);
	S = init; E = init; I = init; R = init; a = init; b = init;
	S(1) = S0; E(1) = E0; I(1) = I0; R(1) = R0;
	Sp = init; Ep = init; Ip = init; Rp = init;
	
	% computation of coefficients a_k and b_k
	for k = 1:N
	    b(k) = k^alpha-(k-1)^alpha;
	    a(k) = (k+1)^(alpha+1)-2*k^(alpha+1)+(k-1)^(alpha+1);
	end
	
	for j = 2:N
	
	    % First part: predict
	
	    % differential equations of control system
	    aux_s = 0; aux_e = 0; aux_i = 0; aux_r = 0;
	    for k = 1:j
	        aux_s = aux_s+b(j-k+1)*(lambda(k)-miu_*S(k)...
	                -beta(k)*S(k)*I(k)+gama_*R(k));
	        aux_e = aux_e+b(j-k+1)*(beta(k)*S(k)*I(k)...
	                -(miu_+epsilon_)*E(k));
	        aux_i = aux_i+b(j-k+1)*(epsilon_*E(k)-(miu_+niu_+Ta(k))*I(k));
	        aux_r = aux_r+b(j-k+1)*(niu_*I(k)-miu_*R(k)-gama_*R(k)...
	                +Ta(k)*I(k));
	    end
	
	    Sp(j) = S0+h^alpha/gamma(1+alpha)*aux_s;
	    Ep(j) = E0+h^alpha/gamma(1+alpha)*aux_e;
	    Ip(j) = I0+h^alpha/gamma(1+alpha)*aux_i;
	    Rp(j) = R0+h^alpha/gamma(1+alpha)*aux_r;
	
	    % Second part: correct
	
	    aux_ss = lambda(j)-miu_*Sp(j)-beta(j)*Sp(j)*Ip(j)+gama_*Rp(j);
	    aux_ee = beta(j)*Sp(j)*Ip(j)-(miu_+epsilon_)*Ep(j);
	    aux_ii = epsilon_*Ep(j)-(miu_+niu_+Ta(j))*Ip(j);
	    aux_rr = niu_*Ip(j)-miu_*Rp(j)-gama_*Rp(j)+Ta(j)*Ip(j);
	
	    auxx = ((j-1)^(alpha+1)-(j-1-alpha)*j^alpha);
	    aux_s0 = auxx*(lambda(1)-miu_*S(1)-beta(1)*S(1)*I(1)+gama_*R(1));
	    aux_e0 = auxx* (beta(1)*S(1)*I(1)-(miu_+epsilon_)*E(1));
	    aux_i0 = auxx*(epsilon_*E(1)-(miu_+niu_+Ta(1))*I(1));
	    aux_r0 = auxx*(niu_*I(1)-miu_*R(1)-gama_*R(1)+Ta(1)*I(1));
	
	    aux_s = 0; aux_e = 0; aux_i = 0; aux_r = 0;
	    for k = 1:j-1
	        aux_s = aux_s+a(j-k)*(lambda(k)-miu_*S(k)-beta(k)*S(k)*I(k)...
	                +gama_*R(k));
	        aux_e = aux_e+a(j-k)*(beta(k)*S(k)*I(k)-(miu_+epsilon_)*E(k));
	        aux_i = aux_i+a(j-k)*(epsilon_*E(k)-(miu_+niu_+Ta(k))*I(k));
	        aux_r = aux_r+a(j-k)*(niu_*I(k)-miu_*R(k)...
	                -gama_*R(k)+Ta(k)*I(k));
	    end
	
	    S(j) = S0+h^alpha/gamma(2+alpha)*(aux_ss+aux_s0+aux_s);
	    E(j) = E0+h^alpha/gamma(2+alpha)*(aux_ee+aux_e0+aux_e);
	    I(j) = I0+h^alpha/gamma(2+alpha)*(aux_ii+aux_i0+aux_i);
	    R(j) = R0+h^alpha/gamma(2+alpha)*(aux_rr+aux_r0+aux_r);
	end
	
	y(1,:) = S; y(2,:) = E; y(3,:) = I; y(4,:) = R;
	
end
	
	
% function III: resolution of the fractional adjoint system
	
	
function [y] = system2_adjoint(S,I,Ta,t,N,alpha)
	
	global  miu gama epsilon niu tfinal k1 b0 b1 phi;
	
	% time-dependent parameter
	fbeta = @(t) b0^alpha.* (1 + b1 * cos( 2 * pi * t + phi ) );
	
	% Correction of values of parameters
	miu_=miu^alpha; niu_=niu^alpha; epsilon_=epsilon^alpha;
	gama_=gama^alpha;
	
	% initialization of variables
	beta = fbeta(t);
	h = tfinal/N; init = zeros(1,N); a = init; b = init;
	p1 = init; p2 = init; p3 = init; p4 = init;
	p1p = init; p2p = init; p3p = init; p4p = init;
	
	% First part: predict
	
	S = S(end:-1:1); I = I(end:-1:1);
	Ta = Ta(end:-1:1); beta = beta(end:-1:1);
	
	% computation of coefficients a_k and b_k
	for k = 1:N
	    b(k) = k^alpha-(k-1)^alpha;
	    a(k) = (k+1)^(alpha+1)-2*k^(alpha+1)+(k-1)^(alpha+1);
	end
	
	for j = 2:N
	
	    % differential equations of adjoint system
	    aux_p1 = 0; aux_p2 = 0; aux_p3 = 0; aux_p4 = 0;
	    for k = 1:j
	        aux_p1 = aux_p1+b(j-k+1)*(-1)*(p1(k)*(miu_+beta(k)*I(k))- ...
	                beta(k)*I(k)*p2(k));
	        aux_p2 = aux_p2+b(j-k+1)*(-1)*(p2(k)*(miu_+epsilon_)...
	                -epsilon_*p3(k));
	        aux_p3 = aux_p3+b(j-k+1)*(-1)*(-k1+beta(k)*p1(k)*S(k)...
	                -p2(k)*beta(k)*S(k)+p3(k)*(miu_+niu_+Ta(k))...
	                -p4(k)*(niu_+Ta(k)));
	        aux_p4 = aux_p4+b(j-k+1)*(-1)*(-gama_*p1(k)...
	                +p4(k)*(miu_+gama_));
	    end
	
	    p1p(j) = h^alpha/gamma(1+alpha)*aux_p1;
	    p2p(j) = h^alpha/gamma(1+alpha)*aux_p2;
	    p3p(j) = h^alpha/gamma(1+alpha)*aux_p3;
	    p4p(j) = h^alpha/gamma(1+alpha)*aux_p4;
	
	    % Second part: correct
	
	    aux_pp1 = (-1)*(p1p(j)*(miu_+beta(j)*I(j))-beta(j)*I(j)*p2p(j));
	    aux_pp2 = (-1)*(p2p(j)*(miu_+epsilon_)-epsilon_*p3p(j));
	    aux_pp3 = (-1)*(-k1+beta(j)*p1p(j)*S(j)-p2p(j)*beta(j)*S(j)...
	                +p3p(j)*(miu_+niu_+Ta(j))-p4p(j)*(niu_+Ta(j)));
	    aux_pp4 = (-1)*(-gama_*p1p(j)+p4p(j)*(miu_+gama_));
	
	    auxx = (-1)*((j-1)^(alpha+1)-(j-1-alpha)*j^alpha);
	    aux_p10 = auxx*(p1(1)*(miu_+beta(1)*I(1))-beta(1)*I(1)*p2(1));
	    aux_p20 = auxx*(p2(1)*(miu_+epsilon_)-epsilon_*p3(1));
	    aux_p30 = auxx*(-k1+beta(1)*p1(1)*S(1)-p2(1)*beta(1)*S(1)...
	                +p3(1)*(miu_+niu_+Ta(1))-p4(1)*(niu_+Ta(1)));
	    aux_p40 = auxx*(-gama_*p1(1)+p4(1)*(miu_+gama_));
	
	    aux_p1 = 0; aux_p2 = 0; aux_p3 = 0; aux_p4 = 0;
	    for k = 1:j-1
	        aux_p1 = aux_p1+a(j-k)*(-1)*(p1(k)*(miu_+beta(k)*I(k))- ...
	                beta(k)*I(k)*p2(k));
	        aux_p2 = aux_p2+a(j-k)*(-1)*( p2(k)*(miu_+epsilon_)...
	                -epsilon_*p3(k));
	        aux_p3 = aux_p3+a(j-k)*(-1)*( -k1+beta(k)*p1(k)*S(k)...
	                -p2(k)*beta(k)*S(k)+p3(k)*(miu_+niu_+Ta(k))...
	                -p4(k)*(niu_+Ta(k)));
	        aux_p4 = aux_p4+a(j-k)*(-1)*(-gama_*p1(k)+p4(k)*(miu_+gama_));
	    end
	
	    p1(j) = h^alpha/gamma(2+alpha)*(aux_pp1+aux_p10+aux_p1);
	    p2(j) = h^alpha/gamma(2+alpha)*(aux_pp2+aux_p20+aux_p2);
	    p3(j) = h^alpha/gamma(2+alpha)*(aux_pp3+aux_p30+aux_p3);
	    p4(j) = h^alpha/gamma(2+alpha)*(aux_pp4+aux_p40+aux_p4);
	end
	
	y(1,:) = p1(end:-1:1); y(2,:) = p2(end:-1:1); y(3,:) = p3(end:-1:1);
	y(4,:) = p4(end:-1:1);
	
end
	
	
% function IV: control projection over the set of admissible controls
	
	
function [v] = projection(vect,trmax)
	
	isNeg = vect<0; vect(isNeg) = 0;
	isHuge = vect>trmax; vect(isHuge) = trmax;
	v = vect;
end
\end{alltt}
\end{appendices}


\begin{adjustwidth}{-\extralength}{0cm}

\reftitle{References}


\PublishersNote{}
\end{adjustwidth}

\end{document}